\documentclass[11pt, twoside]{amsart}

\usepackage[utf8]{inputenc}
\usepackage[T1]{fontenc}
\usepackage{amssymb,amsmath,amstext}
\usepackage{hyperref, graphicx}
\usepackage{stackrel}

\newtheorem{theor}{Theorem}[section]
\newtheorem{lem}[theor]{Lemma}
\newtheorem{defin}[theor]{Definition}

\newtheorem{prop}[theor]{Proposition} 

\newtheorem{notation}[theor]{Notation}
\newtheorem{exam}[theor]{Example}

\newtheorem{cor}[theor]{Corollary}
\newtheorem{rem}[theor]{Remark}

\newtheorem{assump}[theor]{Assumption}

\numberwithin{equation}{section}

\newcommand{\mr}{\mathrm}

\newcommand{\es}{\emptyset}

\newcommand{\mcA}{\mathcal{A}}
\newcommand{\mcB}{\mathcal{B}}

\newcommand{\mff}{\mathfrak{f}}
\newcommand{\mfg}{\mathfrak{g}}

\newcommand{\mbW}{\mathbf{W}}
\newcommand{\mbX}{\mathbf{X}}
\newcommand{\mbY}{\mathbf{Y}}

\newcommand{\mbbP}{\mathbb{P}}
\newcommand{\mbbN}{\mathbb{N}}
\newcommand{\mbbR}{\mathbb{R}}

\newcommand{\msfC}{\mathsf{C}}

\newcommand{\rng}{\mathrm{rng}}

\title[Asymptotic elimination of aggregation functions]
{A general approach to asymptotic elimination of aggregation functions and
generalized quantifiers}

\author{Vera Koponen and
Felix  Weitkämper}

\address{Vera Koponen, Department of Mathematics, Uppsala University, Sweden.}

\email{vera.koponen@math.uu.se}

\address{Felix Weitkämper,
Ludwig-Maximilians-Universität München, Munich, Germany.}

\email{felix.weitkaemper@lmu.de}

\address{Felix Weitkämper, German University of Digital Science, Germany.}

\email{felix.weitkaemper@german-uds.de}

\date{18 November 2025}

\begin{document}

\begin{abstract}
We consider a logic with truth values in the unit interval and which uses aggregation functions instead of quantifiers,
and we describe a general approach to asymptotic elimination of aggregation functions and, indirectly, of
asymptotic elimination of Mostowski style generalized quantifiers, since such can be expressed by using
aggregation functions. The notion of ``local continuity'' of an aggregation function, 
which we make precise in two (related) ways, plays a central role in this approach.
\end{abstract}

\maketitle

\section{Introduction}

\noindent
The use of (generalized) quantifiers is a way to increase the expressivity of a logic, 
beyond the use of functions (often called connectives)
which, for some $k$, assign a truth value to any $k$-tuple of truth values.
In the context of artificial intelligence and machine learning it makes sense to use 
aggregation functions, such as for example the average of a sequence of numbers,
to increase expressivity, for example when
defining probability distributions, or when defining queries. 
(See e.g. \cite[pages 10, 117, 382]{BKNP} or \cite[pages 31, 54]{DKNP}.)
In the context of this article, an aggregation function is a function that takes, for some integer $k >0$,
a $k$-tuple $(\bar{p}_1, \ldots, \bar{p}_k)$ as input, where each $\bar{p}_i$ is a finite sequence of reals
from the unit interval, and gives a real number in the unit interval as output. In addition, for each $i$, the order of
the reals in $\bar{p}_i$ should not matter for the output value. 
(Therefore some people prefer to view each $\bar{p}_i$ as a finite multiset.)
Note that we allow $\bar{p}_i$ and $\bar{p}_j$ to have different length (if $i \neq j$).
To systematize the use of aggregation functions
one can incorporate them into the syntax and semantics of a formal logic, just as done with quantifiers.
But even if a  sequence of numbers contains only the numbers `0' and `1'
the average of the sequence may well be a number strictly between 0 and 1.
So if we use aggregation functions in a formal logic we have to accept 
that the (truth) values of formulas can be other numbers than 0 or 1.

In this work we consider a quite general (formal) logic, which we call $PLA^*$
(``probability logic with aggregation functions''), with truth values in the unit interval $[0, 1]$
and which uses all possible aggregation functions. (The `$*$' in $PLA^*$ is there to indicate that it is another, more
general and flexible, version of the logics $PLA$ and $PLA^+$ considered in \cite{KW1} and \cite{KW2}.)
The restriction of truth values to $[0, 1]$ is partly motivated by the fact that it is often natural to see such truth values
as degrees of uncertainty (e.g. probabilities), 
and partly because a sequence of nonnegative numbers can be ``normalized'' so that the 
resulting sequence contains only numbers in the unit interval but preserves the variation in the original sequence.
Example~\ref{example of page rank} illustrates the expressivity of $PLA^*$ by showing that every ``stage'' 
(in the iterative approximation) of
the PageRank \cite{BP} is definable in $PLA^*$.
The first-order quantifiers $\exists$ and $\forall$ are not lost by using aggregation functions instead of quantifiers,
because $\exists$ and $\forall$ can be expressed by using the aggregation functions `maximum' 
and `minimum' (see Remark~\ref{FO is expressible in PLA*}).
Moreover, as we will see in
Proposition~\ref{generalized quantifiers can be expressed in PLA*},
every condition which can be expressed by first-order logic extended by any generalized quantifiers in
the sense of Mostowski \cite{Mos57, Mos} can be expressed by a formula in $PLA^*$.
(Other early work on generalized quantifiers include Lindström \cite{Lin} and 
Hajek, Havel, and Chytil \cite{HHC} about generating plausible hypotheses in the early days of artificial intelligence.)

We now consider the following context:
$\sigma$ is a finite and relational signature (vocabulary),
$D_n$, $n = 1, 2, 3, \ldots$, is a sequence of finite domains (i.e. sets) such that, for every $n$, 
the cardinality of $D_{n+1}$ is greater than the cardinality of $D_n$,
$\mbW_n$ is a set of $\sigma$-structures with domain $D_n$, and $\mbbP_n$ is
a probability distribution on $\mbW_n$. 
Results have been proven, for various 2-valued logics, various sets of structures $\mbW_n$,
and various sequences of probability distributions $\mbbP_n$, which say that all (or some) formulas of the logic are,
with high probability, equivalent to a formula without quantifiers (or that this is not the case). 
Examples include \cite{DG, FGT, Gleb, Jae98a, Kai, KL, Kop20}.
Often such a result can be used to prove a so-called logical convergence law, or even a zero-one law.
This kind of result often has implications for the expressivity of the logic considered and for
the existence of computationally efficient ways of estimating, on a large domain, the probability of a query defined by the logic.

Within the specified context we can ask, for $\varphi(\bar{x}) \in PLA^*$, $0 \leq c < d \leq 1$, and $\bar{a} \in (D_n)^{|\bar{x}|}$,
what the probability is that the value of $\varphi(\bar{a})$ in a random $\mcA \in \mbW_n$ is in the interval $(c, d)$.
Let us call a formula of any logic {\em aggregation-free} if it contains no aggregation function and no quantifier.
For the reasons mentioned above, it is of interest to understand
under what conditions for a $PLA^*$-formula $\varphi(\bar{x})$ there is a $PLA^*$-formula $\psi(\bar{x})$ such that $\psi(\bar{x})$ is 
aggregation-free, or at least of smaller complexity than $\varphi$, and for all $\varepsilon > 0$
the probability that the values of the two formulas differ by at most $\varepsilon$ (for all substitutions of parameters
for the free variables $\bar{x}$) tends to 1 as $n\to\infty$; in this case we say that $\varphi(\bar{x})$ 
and $\psi(\bar{x})$ are asymptotically equivalent
(see Definition~\ref{definition of asymptotic equivalence}).

Such results were proved in \cite{KW1, KW2}
and the proofs in both articles have a common strategy.
This strategy also appears to some extent in some of the proofs in \cite{Jae98a, KL, Kop20},
although it is not as clear because the later articles deal with quantifiers instead of aggregation functions.
In a simplified description, the first step of the strategy is to prove that for some set $L_0$
of  0/1-valued aggregation-free formulas, 
if $\varphi(\bar{x}, \bar{y}) \in L_0$ then there are finitely many $\varphi_i(\bar{x}) \in L_0$ such that the following holds
with probability tending to 1 as $n\to\infty$:
If $\bar{a}$ is a sequence of elements from $D_n$ of the same length as $\bar{x}$, 
then $\varphi_i(\bar{a})$ holds for some $i$, and in this case
the proportion of $\bar{b}$ (among tuples of elements from $D_n$ with the same length as $\bar{y}$)
such that $\varphi(\bar{a}, \bar{b})$ holds is, with high probability, close to some number $\alpha_i$ depending only on
$\varphi$, $\varphi_i$ and the sequence of probability distributions $\mbbP_n$.
The second step is to use the conclusions of the first step to show that for every formula $\varphi(\bar{x})$
there is a formula $\psi(\bar{x})$ without aggregation functions and without quantifiers which is 
asymptotically equivalent to $\varphi(\bar{x})$.

In general, the first step seems to be intricately connected to the kind of structures and the kind of probability distributions considered.
But the second step can be formulated
in such a way that it makes sense in a wide variety of situations,  as we will show in 
Assumption~\ref{assumptions on the basic logic} and 
Theorem~\ref{main result}.
We hope that Theorem~\ref{main result} can be used several times in the future once its preconditions 
(Assumption~\ref{assumptions on the basic logic}) have been
shown to hold. 
This avoids carrying out a similar proof again and again, 
as already done in \cite{KW1} and \cite{KW2}, and also in e.g. \cite{KL} and \cite{Kop20}, 
although it is less clear in the latter cases. 
In fact, Theorem~\ref{main result} has been used in precisely this way in \cite{Kop26} and \cite{KT}
(which were written after the submission of this article).

Theorem~\ref{main result} also gives us a general understanding of
conditions under which aggregation functions can
be asymptotically eliminated from $PLA^*$-formulas.
Indirectly it also gives us an understanding of when Mostowski style generalized quantifiers can be asymptoticaly eliminated. Hence our Theorem~\ref{main result} is similar in spirit, but different in the technical approach, to some results of Kaila \cite[Theorems~4.4 and~4.7]{Kai}.

In some contexts it may be impossible to find, for a formula $\varphi(\bar{x}) \in PLA^*$
an aggregation-free formula $\psi(\bar{x}) \in PLA^*$ which is asymptotically equivalent to $\varphi(\bar{x})$.
Theorem~\ref{main result} may still be useful and show that $\varphi(\bar{x})$ is asymptotically equivalent
to a formula $\psi(\bar{x})$ which is, in some sense, simpler than $\varphi(\bar{x})$.
For example, $\psi(\bar{x})$ may be a formula with only ``bounded aggregation functions (or quantifiers)''
when this notion makes sense (as e.g. in the situation where some relations have bounded Gaifman degree).

The discussion so far has hidden the fact that, in addition to certain preconditions
stated in Assumption~\ref{assumptions on the basic logic}, a sufficient and (in general) necessary condition for asymptotically
eliminating an aggregation function $F$ in a formula $\varphi \in PLA^*$ is that $F$ is, in a sense, {\em continuous}
on a ``local'' set which is
determined by some parameters that are determined by the subformulas of $\varphi$ and the 
sequence of probability distributions $\mbbP_n$.

Section~\ref{Notions of continuity for aggregation functions}
focuses on two continuity properties of aggregation functions that are crucial to this study.
The first continuity property has the advantage that it is relatively easy to check whether an aggregation function has it.
The other property is formulated so that it is easy to use in the proof of Theorem~\ref{main result}.
However, as stated by Proposition~\ref{proposition relating the notions of continuity},
both ``local'' continuity properties are very closely related, possibly equivalent although we have not been able to show this.
The ``global'' versions of the two continuity properties are indeed equivalent.

The structure of the article is as follows:
Section~\ref{Logic with aggregation functions} defines the notions of connective, aggregation function and the
logic $PLA^*$ that we will work with.
Section~\ref{Expressing quantifiers with aggregation functions} defines generalized quantifiers in the sense
of Mostowski and shows that they can be represented by aggregation functions.
Section~\ref{Notions of continuity for aggregation functions} defines the two notions of ``local'' continuity
that we consider, shows how they are related, and gives examples.
Finally, Section~\ref{Asymptotic equivalence to basic formulas}
describes a general approach to asymptotic elimination of aggregation functions which is concluded by
Theorem~\ref{main result}.

\subsection*{Notation and terminology}

By $\mbbN$ and $\mbbN^+$ we denote the sets of nonnegative, respectively, positive integers.
We denote finite sequences by $\bar{a}, \bar{b}, \ldots$, $\bar{x}, \bar{y}, \ldots$,
where typically $\bar{x}, \bar{y}, \ldots$ denote finite sequences of distinct variables.
The length of a sequence $\bar{a}$ is denoted by $|\bar{a}|$, and 
the set of elements occuring in $\bar{a}$ is denoted by $\rng(\bar{a})$.
If $A$ is a set then $|A|$ denotes its cardinality, and we let
$A^{<\omega}$ denote the set of all finite nonempty sequences of elements from $A$,
so $A^{<\omega} = \bigcup_{n\in\mbbN^+}A^n$.
First-order structures are denoted by $\mcA, \mcB, \ldots$ and their domains (or universes) by $A, B, \ldots$
unless we say something else (for example, in Section~\ref{Asymptotic equivalence to basic formulas} we consider
a fixed sequence of domains $D_n$, $n \in \mbbN^+$).
We use the word {\em signature} (or vocabulary) in its usual first-order sense.
If a signature contains only relation symbols we call it {\em relational}.

\section{Logic with aggregation functions}\label{Logic with aggregation functions}

\noindent
In this section we define the logic that we will work with, called $PLA^*$.
It is an extension of the logics $PLA$ and $PLA^+$ considered in \cite{KW1} and \cite{KW2}, respectively,
by allowing a more flexible use of aggregation functions. Formulas of $PLA^*$ may take any (truth) value in the
unit interval $[0, 1]$ and we think of the values 0 and 1 as corresponding to ``false'' and ``true''.
$PLA^*$ uses aggregation functions instead of quantifiers to ``aggregate'' over a domain.
In 
Section~\ref{Expressing quantifiers with aggregation functions}
we will see that any condition expressed by a formula of first-order logic extended by
generalized quantifiers in the sense of Mostowski \cite{Mos} can be
expressed by a formula of $PLA^*$.
Recall that $[0, 1]^{<\omega}$ denotes the set of all finite nonempty sequences of reals in 
the unit interval $[0, 1]$.

\begin{defin}\label{definition of connective and aggregation function} {\rm
Let $k \in \mbbN^+$.\\
(i) A function $\msfC : [0, 1]^k \to [0, 1]$ will also be called a {\em ($k$-ary) connective}.\\
(ii) A function $F : \big([0, 1]^{<\omega}\big)^k \to [0, 1]$ which is symmetric in the following sense
will be called a {\em ($k$-ary) aggregation function}:
if $\bar{p}_1, \ldots, \bar{p}_k \in [0, 1]^{<\omega}$ and, for $i = 1, \ldots, k$,
$\bar{q}_i$ is a reordering of the entries in $\bar{p}_i$,
then $F(\bar{q}_1, \ldots, \bar{q}_k) = F(\bar{p}_1, \ldots, \bar{p}_k)$.
}\end{defin}

\noindent
The functions defined in the next definition are continuous and when restricted to $\{0, 1\}$
(as opposed to the interval $[0, 1]$) they have the usual meanings  of $\neg$, $\wedge$, $\vee$, and
$\rightarrow$. (The definitions correspond to the semantics of Lukasiewicz logic 
(see for example \cite[Section~11.2]{Ber}, or \cite{LT}).

\begin{defin}\label{special connectives} {\bf (Some special continuous connectives)} {\rm
Let
\begin{enumerate}
\item $\neg : [0, 1] \to [0, 1]$ be defined by $\neg(x) = 1 - x$,
\item $\wedge : [0, 1]^2 \to [0, 1]$ be defined by $\wedge(x, y) = \min(x, y)$,
\item $\vee : [0, 1]^2 \to [0, 1]$ be defined by $\vee(x, y) = \max(x, y)$, and
\item $\rightarrow : [0, 1]^2 \to [0, 1]$ be defined by $\rightarrow(x, y) = \min(1, \ 1 - x + y)$.
\end{enumerate}
}\end{defin}

\begin{defin}\label{examples of aggregation functions} {\bf (Some common aggregation functions)} {\rm
For $\bar{p} = (p_1, \ldots, p_n) \in [0, 1]^{<\omega}$, define
\begin{enumerate}
\item $\max(\bar{p})$ to be the {\em maximum} of all $p_i$,
\item $\min(\bar{p})$ to be the {\em minimum} of all $p_i$,
\item $\mr{am}(\bar{p}) = (p_1 + \ldots + p_n)/n$, so `am' is the {\em arithmetic mean}, or {\em average}, and
\item $\mr{gm}(\bar{p}) = \big(\prod_{i=1}^n p_i\big)^{(1/n)}$, so `gm' is the {\em geometric mean}.
\end{enumerate}
}\end{defin}

\noindent
In the examples that will follow we will meet more aggregation functions.
{\em For the rest of this section we fix a finite and relational signature $\sigma$.}

\begin{defin}\label{syntax of PLA*}{\bf (Syntax of $PLA^*(\sigma)$)} {\rm 
We define formulas of $PLA^*(\sigma)$, as well as the set of free variables  of a formula $\varphi$, 
denoted $Fv(\varphi)$, as follows.
\begin{enumerate}
\item  For each $c \in [0, 1]$, $c \in PLA^*(\sigma)$ (i.e. $c$ is a formula) and $Fv(c) = \es$. 
We also let $\bot$ and $\top$
denote $0$ and $1$, respectively.

\item For all variables $x$ and $y$, `$x = y$' belongs to $PLA^*(\sigma)$ and $Fv(x = y) = \{x, y\}$.

\item For every $R \in \sigma$, say of arity $r$, and any choice of variables $x_1, \ldots, x_r$, $R(x_1, \ldots, x_r)$ belongs to 
$PLA^*(\sigma)$ and  $Fv(R(x_1, \ldots, x_r)) = \{x_1, \ldots, x_r\}$.

\item If $n \in \mbbN^+$, $\varphi_1, \ldots, \varphi_n \in PLA^*(\sigma)$ and
$\msfC : [0, 1]^n \to [0, 1]$ is a continuous connective, then 
$\msfC(\varphi_1, \ldots, \varphi_n)$ is a formula of $PLA^*(\sigma)$ and
its set of free variables is $Fv(\varphi_1) \cup \ldots \cup Fv(\varphi_n)$.

\item Suppose that $\varphi_1, \ldots, \varphi_k \in PLA^*(\sigma)$,
$\chi_1, \ldots, \chi_k \in PLA^*(\sigma)$,
$\bar{y}$ is a sequence of distinct variables,
and that $F : \big( [0, 1]^{<\omega} \big)^k \to [0, 1]$ is an aggregation function.
Then 
\[
F(\varphi_1, \ldots, \varphi_k : \bar{y} : 
\chi_1, \ldots, \chi_k)
\]
is a formula of $PLA^*(\sigma)$ and its set of free variables is
\[
\big( \bigcup_{i=1}^k Fv(\varphi_i)\big) \setminus \rng(\bar{y}),
\] 
so this construction binds the variables in $\bar{y}$.
\end{enumerate}
}\end{defin}

\begin{notation}\label{abbreviation when using aggregation functions}{\rm
(i) When denoting a formula in $PLA^*(\sigma)$ by e.g. $\varphi(\bar{x})$ we assume that $\bar{x}$ is a sequence
of different variables and that every variable in the formula denoted by $\varphi(\bar{x})$ belongs to $\rng(\bar{x})$
(but we do not require that every variable in $\rng(\bar{x})$ actually occurs in the formula).\\
(ii) If all $\chi_1, \ldots, \chi_k$ are the same formula $\chi$,
then we may abbreviate 
\[
F(\varphi_1, \ldots, \varphi_k : \bar{y} : 
\chi_1, \ldots, \chi_k)
\qquad \text{ by } \qquad
F(\varphi_1, \ldots, \varphi_k : \bar{y} : \chi).
\]
}\end{notation}

\begin{defin}\label{definition of literal} {\rm
The $PLA^*(\sigma)$-formulas described in parts~(2) and~(3) of 
Definition~\ref{syntax of PLA*}
are called {\em first-order atomic formulas}.
A $PLA^*(\sigma)$-formula is called a {\em first-order literal} if it has the form $\varphi(\bar{x})$
or $\neg\varphi(\bar{x})$, where $\varphi(\bar{x})$ is a first-order atomic formula and
$\neg$ is like in Definition~\ref{special connectives}
(so it corresponds to negation when truth values are restricted to 0 and~1).
}\end{defin}

\begin{defin}\label{semantics of PLA*}{\bf (Semantics of $PLA^*(\sigma)$)} {\rm
For every $\varphi \in PLA^*(\sigma)$ and every sequence of distinct variables $\bar{x}$ such that 
$Fv(\varphi) \subseteq \rng(\bar{x})$ we associate a mapping from pairs $(\mcA, \bar{a})$,
where $\mcA$ is a finite $\sigma$-structure and $\bar{a} \in A^{|\bar{x}|}$, to $[0, 1]$.
The number in $[0, 1]$ to which $(\mcA,\bar{a})$ is mapped is denoted by $\mcA(\varphi(\bar{a}))$
and is defined by induction on the complexity of formulas.
\begin{enumerate}
\item If $\varphi(\bar{x})$ is a constant $c$ from $[0, 1]$, then $\mcA(\varphi(\bar{a})) = c$.

\item If $\varphi(\bar{x})$ has the form $x_i = x_j$, then $\mcA(\varphi(\bar{a})) = 1$ if $a_i = a_j$,
and otherwise $\mcA(\varphi(\bar{a})) = 0$.

\item For every $R \in \sigma$, of arity $r$ say, if $\varphi(\bar{x})$ has the form $R(x_{i_1}, \ldots, x_{i_r})$,
then $\mcA(\varphi(\bar{a})) = 1$ if $\mcA \models R(a_{i_1}, \ldots, a_{i_r})$
(where `$\models$' has the usual meaning
of first-order logic), and otherwise $\mcA(\varphi(\bar{a})) = 0$.

\item If $\varphi(\bar{x})$ has the form $\msfC(\varphi_1(\bar{x}), \ldots, \varphi_k(\bar{x}))$,
where $\msfC : [0, 1]^k \to [0, 1]$ is a continuous connective, then
\[
\mcA\big(\varphi(\bar{a})\big) \ = \ 
\msfC\big(\mcA(\varphi_1(\bar{a})), \ldots, \mcA(\varphi_k(\bar{a}))\big).
\]

\item Suppose that $\varphi(\bar{x})$ has the form 
\[
F(\varphi_1(\bar{x}, \bar{y}), \ldots, \varphi_k(\bar{x}, \bar{y}) : \bar{y} : 
\chi_1(\bar{x}, \bar{y}), \ldots, \chi_k(\bar{x}, \bar{y}))
\] 
where $\bar{x}$ and $\bar{y}$ are sequences of distinct variables, $\rng(\bar{x}) \cap \rng(\bar{y}) = \es$, and
$F : \big( [0, 1]^{<\omega} \big)^k \to [0, 1]$ is an aggregation function.
If, for every $i = 1, \ldots, k$, the set 
$\{\bar{b} \in A^{|\bar{y}|} : \mcA(\chi_i(\bar{a}, \bar{b})) = 1\}$ is nonempty, then
let 
\[
\bar{p}_i = 
\big(\mcA\big(\varphi_i(\bar{a}, \bar{b})\big) : \bar{b} \in A^{|\bar{y}|} \text{ and } 
\mcA\big(\chi_i(\bar{a}, \bar{b})\big) = 1\big)
\]
and 
\[
\mcA\big(\varphi(\bar{a})\big) = F(\bar{p}_1, \ldots, \bar{p}_k).
\]
Otherwise let $\mcA\big(\varphi(\bar{a})\big) = 0$.
\end{enumerate}
}\end{defin}

\begin{defin}\label{definition of set defined by a formula} {\rm
(i) Suppose that $\varphi(\bar{x}, \bar{y}) \in PLA^*(\sigma)$,
$\bar{x}$ and $\bar{y}$ are sequences of distinct variables, $\rng(\bar{x}) \cap \rng(\bar{y}) = \es$
$\mcA$ is a finite $\sigma$-structure and $\bar{a} \in A^{|\bar{x}|}$.
Then $\varphi(\bar{a}, \mcA)$ denotes the set 
$\big\{\bar{b} \in A^{|\bar{y}|} : \mcA\big(\varphi(\bar{a}, \bar{b})\big) = 1\big\}$.\\
(ii) If $\varphi(\bar{x}) \in PLA^*(\sigma)$, $\mcA$ is a finite $\sigma$-structure and that $\bar{a} \in A^{|\bar{x}|}$,
then `$\mcA \models \varphi(\bar{a})$' means the same as `$\mcA\big(\varphi(\bar{a})\big) = 1$'.
}\end{defin}

\begin{exam}\label{example of page rank} {\rm
As an example of what can be expressed with $PLA^*(\sigma)$ we consider the notion of PageRank \cite{BP}.
The PageRank of an internet site can be approximated in ``stages'' as follows
(if we supress the ``damping factor'' for simplicity), where $IN_x$ is the set of sites that link to $x$,
and $OUT_y$ is the set of sites that $y$ link to:
\begin{align*}
&PR_0(x) = 1/N \text{ where $N$ is the number of sites,} \\
&PR_{k+1}(x) = \sum_{y \in IN_x} \frac{PR_k(y)}{|OUT_y|}.
\end{align*}
It is not difficult to prove, by induction on $k$, that for every $k$ the sum of all $PR_k(x)$ as $x$ ranges over all sites
is 1. Hence the sum in the definition of $PR_{k+1}$ is less or equal to 1 (which will matter later).
Let $E \in \sigma$ be a binary relation symbol representing a link.
Define the aggregation function $\mr{length}^{-1} : [0, 1]^{<\omega} \to [0, 1]$ by $\mr{length}^{-1}(\bar{p}) = 1/|\bar{p}|$.
Then $PR_0(x)$ is expressed by the $PLA^*(\sigma)$-formula $\mr{length}^{-1}(x = x : y : \top)$.

Suppose that $PR_k(x)$ is expressed by $\varphi_k(x) \in PLA^*(\sigma)$.
Note that multiplication is a continuous connective from $[0, 1]^2$ to $[0, 1]$ so it can be used in $PLA^*(\sigma)$-formulas.
Then observe that the quantity $|OUT_y|^{-1}$ is expressed by the $PLA^*(\sigma)$-formula
\[
\mr{length}^{-1}\big(y=y : z : E(y, z)\big)
\]
which we denote by $\psi(y)$.
Let $\mr{tsum} : [0, 1]^{<\omega} \to [0, 1]$ be the ``truncated sum'' defined by letting
$\mr{tsum}(\bar{p})$ be the sum of all entries in $\bar{p}$ if the sum is at most 1, and otherwise $\mr{tsum}(\bar{p}) = 1$.
Then $PR_{k+1}(x)$ is expressed by the $PLA^*(\sigma)$-formula
\[
\mr{tsum}\big(x = x \wedge (\varphi_k(y) \cdot \psi(y)) : y : E(y, x)\big).
\]
With $PLA^*$ we can also define all stages of the SimRank \cite{JW} in a simpler way than done
in \cite{KW1} with the sublogic $PLA \subseteq PLA^*$.
}\end{exam}

\begin{rem}\label{FO is expressible in PLA*} {\bf (The relation to first-order logic)} {\rm
From the syntax and semantics of $PLA^*(\sigma)$ and the fact
the connectives $\neg$, $\wedge$, $\vee$ and $\to$ are continuous
it follows that every quantifier-free first-order formula
$\varphi(\bar{x})$ (over $\sigma$) is also a $PLA^*(\sigma)$-formula.
Now suppose that $\varphi(\bar{x}, y)$ is a first-order formula and $\psi(\bar{x}, y)$ a $PLA^*(\sigma)$-formula 
such that for every finite $\sigma$-structure $\mcA$ and $\bar{a} \in A^{|\bar{x}|}$
and $b \in A$,
$\mcA \models \varphi(\bar{a}, b)$ if and only if $\mcA\big(\psi(\bar{a}, b)\big) = 1$.
Then, for all $\bar{a} \in A^{|\bar{x}|}$,
\[
\mcA \models \exists y \varphi(\bar{a}, y) \text{ if and only if } 
\mcA\big(\max(\psi(\bar{a}, y) : y : \top)\big) = 1.
\]
Similarly, the quantifier $\forall$ can be expressed in $PLA^*$ by using  the aggregation function min.
By induction on the complexity of first-order formulas it follows that, for relational signatures and finite structures,
every condition that can be expressed by first-order logic can be expressed with $PLA^*$.
}\end{rem}

\begin{notation}\label{using first-order notation} {\bf (Using $\exists$ and $\forall$ as abbreviations)} {\rm
Motivated by the discussion in Remark~\ref{FO is expressible in PLA*}, if $\varphi(\bar{x}, \bar{y}) \in PLA^*(\sigma)$ 
is a formula which can only take the values 0 or 1 (e.g a boolean combination of atomic first-order formulas)
then we may
use `$\exists \bar{y} \varphi(\bar{x}, \bar{y})$' to mean the same as
`$\max(\varphi(\bar{x}, \bar{y}) : \bar{y} : \top)$', and
`$\forall \bar{y} \varphi(\bar{x}, \bar{y})$' to mean the same as
`$\min(\varphi(\bar{x}, \bar{y}) : \bar{y} : \top)$'.
}\end{notation}

\section{Expressing generalized quantifiers with aggregation functions}
\label{Expressing quantifiers with aggregation functions}

\noindent
In this section we consider generalized quantifiers in the sense of Mostowski \cite{Mos} and
see that they can be expressed by aggregation functions. 
Let $\sigma$ be a relational signature.

\begin{defin}\label{definition of generalized quantifier}{\rm
A {\em quantifier aggregating $k$ sets} is a class $Q$ consisting of 
tuples $(D, X_1, \ldots, X_k)$ such that $D$ is a set, $X_1, \ldots, X_k \subseteq D$ and the following condition holds:
\begin{itemize}
\item[] If $|D| = |E|$, $X_i \subseteq D$, $Y_i \subseteq E$, $|X_i| = |Y_i|$, for $i = 1, \ldots, k$, and
$(D, X_1, \ldots, X_k) \in Q$, then $(E, Y_1, \ldots, Y_k) \in Q$.
\end{itemize}
}\end{defin}

\begin{defin}\label{definition of logic with generalized quantifiers} {\rm
(i) Then let $FOGQ(\sigma)$ be the set of all expressions that can be obtained by adding
the following construction to the definition of a first-order formula over $\sigma$:
\begin{enumerate}
\item[] If $\varphi_1(\bar{x}, \bar{y}), \ldots, \varphi_k(\bar{x}, \bar{y}) \in FOGQ(\sigma)$,
$\bar{y}$ is a sequence of different variables, and 
$Q$ is a quantifier aggregating $k$ sets, then the following expression belongs to $FOGQ(\sigma)$:
\begin{equation*}\label{expression to be interpreted}
Q\bar{y}(\varphi_1(\bar{x}, \bar{y}), \ldots, \varphi_k(\bar{x}, \bar{y})).
\end{equation*}
\end{enumerate}
(ii) The semantics of first-order logic is now extended to $FOGQ(\sigma)$ as follows if $|\bar{y}| = n$:
\begin{enumerate}
\item[] $\mcA \models Q\bar{y} (\varphi_1(\bar{a}, \bar{y}), \ldots, \varphi_k(\bar{a}, \bar{y}))\big)$
if and only if
$(A^n, X_1, \ldots, X_k) \in Q$ where, for all $i = 1, \ldots, k$,
$X_i = \{\bar{b} \in A^n : \mcA \models \varphi_i(\bar{a}, \bar{b})) \}$.
\end{enumerate}
}\end{defin}

\begin{prop}\label{generalized quantifiers can be expressed in PLA*}
Let $\varphi(\bar{x}) \in FOGQ(\sigma)$.
Then there is $\psi(\bar{x}) \in PLA^*(\sigma)$ such that for every finite $\sigma$-structure
$\mcA$ and every $\bar{a} \in A^{|\bar{a}|}$,
$\mcA \models \varphi(\bar{a})$ if and only if $\mcA\big(\psi(\bar{a})\big) = 1$,
and $\mcA \not\models \varphi(\bar{a})$ if and only if $\mcA\big(\psi(\bar{a})\big) = 0$.
\end{prop}

\noindent
{\bf Proof.}
For any $\varphi(\bar{x}) \in FOGQ(\sigma)$, finite $\sigma$-structure $\mcA$ and $\bar{a} \in A^{|\bar{x}|}$,
let `$\mcA\big(\varphi(\bar{a})\big) = 1$' mean the same as `$\mcA \models \varphi(\bar{a})$' and let
'$\mcA\big(\varphi(\bar{a})\big) = 0$' mean the same as `$\mcA \not\models \varphi(\bar{a})$'.

In Remark~\ref{FO is expressible in PLA*} 
we saw that the proposition holds if $\varphi(\bar{x})$ is a first-order formula.
In the same remark we also saw that $\neg$, $\wedge$, $\vee$ and $\to$ can, with Lukasiewicz semantics, be expressed by 
continuous functions from $[0, 1]$ or $[0, 1]^2$ to $[0, 1]$, and $\exists$ and $\forall$ can be expressed by
the aggregation functions max and min, respectively.

Hence is suffices to show the following:

\medskip
\noindent
{\bf Claim:}
{\em Suppose that $\bar{x}$ and $\bar{y}$ are sequences of different variables,
$\rng(\bar{x}) \cap \rng(\bar{y}) = \es$, 
$\varphi_1(\bar{x}, \bar{y}),  \ldots, \varphi_k(\bar{x}, \bar{y}) \in FOGQ(\sigma)$,
and $\psi_1(\bar{x}, \bar{y}), \ldots, \psi_k(\bar{x}, \bar{y}) \in PLA^*(\sigma)$ are such that
for every finite $\sigma$-structure $\mcA$, all $\bar{a} \in A^{|\bar{x}|}$, and all $\bar{b} \in A^{|\bar{y}|}$,
\[
\mcA\big(\varphi_i(\bar{a}, \bar{b})\big) = \mcA\big(\psi_i(\bar{a}, \bar{b})\big) \text{ for all $i = 1, \ldots, k$.}
\]
If $Q$ is an $n$-ary quantifier aggregating $k$ sets, then there is an aggregation function
$F : \big([0, 1]^{<\omega}\big)^k \to [0, 1]$  such that
for every finite $\sigma$-structure $\mcA$ and $\bar{a} \in A^{|x|}$:
\begin{equation}\label{same truth value for quantifier and aggregation function}
\mcA\big(Q\bar{y} (\varphi_1(\bar{a}, \bar{y}), \ldots, \varphi_k(\bar{a}, \bar{y}))\big) \ = \
\mcA\big(F\big(\psi_1(\bar{a}, \bar{y}), \ldots, \psi_k(\bar{a}, \bar{y}) : \bar{y} : \top \big)\big), 
\end{equation}
where we recall that the $PLA^*$-formula `$\top$' (or `1') has the value 1 in every structure.
}

\medskip

\noindent
Let $F : \big([0, 1]^{<\omega}\big)^k \to [0, 1]$  be defined as follows:
If $\bar{p}_i = (p_{i, 1}, \ldots, p_{i, m_i})$ for $i = 1, \ldots, k$, then
\begin{align*}
&F(\bar{p}_1, \ldots, \bar{p}_k) = \\
&\begin{cases}
1 \qquad \text{ if } (D, X_1, \ldots, X_k) \in Q 
\text{ where $D = [m]$, $m = \max(m_1, \ldots, m_k)$} \\ 
\qquad \quad \text{ and } X_i = \{j \in [m_i] : r_{i, j} = 1\}, \text{ and} \\
0 \qquad \text{ otherwise.}
\end{cases}
\end{align*}
For every finite structure $\mcA$ and every $\bar{a} \in A^{|\bar{x}|}$ we now have
\begin{align*}
&\mcA\big(Q\bar{y}(\varphi_1(\bar{a}, \bar{y}), \ldots, \varphi_k(\bar{a}, \bar{y}))\big) = 1 \Longleftrightarrow \\
&\mcA \models Q\bar{y}(\varphi_1(\bar{a}, \bar{y}), \ldots, \varphi_k(\bar{a}, \bar{y})) \Longleftrightarrow \\
&(A^{|\bar{y}|}, \varphi_1(\bar{a}, \mcA), \ldots, \varphi_k(\bar{a}, \mcA)) \in Q \Longleftrightarrow \\
&F(\bar{p}_1, \ldots, \bar{p}_k) = 1
\text{ where, for $i = 1, \ldots, k$, 
$\bar{p}_i = \big(\mcA\big(\varphi_i(\bar{a}, \bar{b})\big) : \bar{b} \in A^{|\bar{y}|}\big)$.}
\end{align*}
Since also
\[
\mcA\big(F\big(\varphi_1(\bar{a}, \bar{y}), \ldots, \varphi_k(\bar{a}, \bar{y}) : \bar{y} : \top \big)\big) = 
F(\bar{p}_1, \ldots, \bar{p}_k)
\]
the equation~(\ref{same truth value for quantifier and aggregation function}) follows.
This proves the claim and concludes the proof of the proposition.
\hfill $\square$

\section{Notions of continuity for aggregation functions}\label{Notions of continuity for aggregation functions}

\noindent
Our goal is to formulate conditions that apply to a wide variety of contexts and under which 
an aggregation function in a $PLA^*$-formula can be ``asymptotically eliminated''.
If all aggregation functions in the formula can be asymptotically eliminated then we get
a ``simpler'' formula which is ``asymptotically equivalent''
(Definition~\ref{definition of asymptotic equivalence})
to the original formula.
One of the conditions that need (necessarily according to 
Remark~\ref{necessity of continuity}) 
to be satisfied is that the aggregation function to be asymptotically eliminated,
say $F : \big([0, 1]^{<\omega}\big)^m \to [0, 1]$,
has some sort of continuity property on a subset of $\big([0, 1]^{<\omega}\big)^m$.

We formulate two notions of continuity, {\em ct-continuity} and {\em up-continuity}, which are ``local'' versions
of the ``global'' notions of  {\em (strongly) admissible aggregation function}, respectively,
{\em (strongly) admissible aggregation function sensu novo} that are considered in \cite{KW1, KW2}.
Besides being localizations of the corresponding notions in \cite{KW1, KW2}, the notions considered here 
are a little bit weaker (than the corresponding global notions in \cite{KW1, KW2}) also in the sense that 
some small conditions in the definitions
of (strong) admissibility and (strong) admissibility sensu novo are 
not present in the corresponding localizations considered here, because
we have realized that those conditions are not necessary for proving ``asymptotic elimination results''.

Then we show
(see Proposition~\ref{proposition relating the notions of continuity})
that ct-continuity and up-continuity are closely related. 
Nevertheless we consider both notions useful because ct-continuity (with respect to some parameters) is
usually easier to verify for concrete aggregation functions, while up-continuity is tailored for making the
proof of ``elimination of an aggregation function'' (Theorem~\ref{main result}) work out.

Throughout the section we give examples of aggregation functions that are ct-continuous, respectively up-continuous,
with respect to some (or all) parameters.

\subsection{Convergence test continuity}

\noindent
In this section we define the local version of {\em (strong) admissibility} (used in \cite{KW1, KW2}) which we call
{\em ct-continuity with respect to a sequence of parameters}.
Before doing that we must define the notion of {\em convergence testing sequence} which generalizes a
notion used by Jaeger in \cite{Jae98a}.
The intuition is that a sequence $\bar{p}_n \in [0, 1]^{<\omega}$, $n \in \mbbN^+$ is convergence testing for
for parameters 
$c_1, \ldots, c_k \in [0,1]$ and $\alpha_1, \ldots \alpha_k  \in  [0,1]$
if the length of $\bar{p}_n$ tends to infinity as $n\to\infty$ and, as $n\to\infty$,
all entries in $\bar{p}_n$ congregate ever closer to the ``convergence points" in the set $\{c_1, \ldots, c_k\}$,
and the proportion of entries in $\bar{p}$ which are close to $c_i$ converges to $\alpha_i$.

\begin{defin}\label{definition of convergence testing} {\bf (Convergence testing sequence)} {\rm 
A sequence $\bar{p}_n \in [0, 1]^{<\omega}$, $n \in \mbbN$,  is called {\em convergence testing} for parameters 
$c_1, \ldots, c_k \in [0,1]$ and $\alpha_1, \ldots \alpha_k  \in  [0,1]$ if the following hold, 
where $p_{n,i}$ denotes the $i$th entry of $\bar{p}_n$:
\begin{enumerate}
\item $|\bar{p}_n| < |\bar{p}_{n+1}|$ for all $n \in \mbbN$.
\item For every disjoint family of open (with respect to the induced topology on $[0, 1]$)
intervals $I_1, \ldots I_k \subseteq [0,1]$ 
 such that $c_i \in I_i$ for each $i$, 
there is an $ N \in \mbbN$ such that $\mathrm{rng}(\bar{p}_n) \subseteq \bigcup\limits_{j=1}^{k} I_j$ for all $n \geq N$, 
and for every $j \in \{1, \ldots, k \}$, 
\[
\lim\limits_{n \rightarrow \infty} \frac{\left| \{ i  : p_{n,i} \in I_j \} \right| }{|\bar{p}_n|} = \alpha_j.
\]
\end{enumerate}   

More generally, a sequence of $m$-tuples 
$(\bar{p}_{1, n}, \ldots, \bar{p}_{m, n}) \in  \big([0, 1]^{<\omega}\big)^m$, $n \in \mbbN$,  
is called {\em convergence testing} for parameters $c_{i,j} \in [0,1]$ and $\alpha_{i,j} \in [0,1]$, 
where $i \in \{1, \ldots, m\}$, $j \in \{ 1, \ldots, k_i \}$ and $k_1, \ldots k_m \in \mbbN^+$, if  for every fixed 
$i \in \{1, \ldots, m \}$
the sequence $\bar{p}_{i, n}$, $n \in \mbbN$, is convergence testing for $c_{i,1}, \ldots, c_{i, k_i}$, and 
$\alpha_{i,1}, \ldots, \alpha_{i, k_i}$.
}\end{defin}

\begin{defin} \label{definition of ct-continuous}{\bf (Convergence test continuity)} {\rm 
An aggregation function \\ 
$F : \big([0, 1]^{<\omega}\big)^m \to [0, 1]$ is called 
{\em ct-continuous (convergence test continuous)} with respect to the sequence of parameters
$c_{i,j}, \alpha_{i,j} \in [0, 1]$, $i = 1, \ldots, m$, $j = 1, \ldots, k_i$, 
if the following condition holds:
\begin{enumerate}
\item[] For all convergence testing sequences of $m$-tuples
$(\bar{p}_{1, n}, \ldots, \bar{p}_{m, n}) \in  \big([0, 1]^{<\omega}\big)^m$, $n \in \mbbN$,
and $(\bar{q}_{1, n}, \ldots, \bar{q}_{m, n}) \in  \big([0, 1]^{<\omega}\big)^m$, $n \in \mbbN$,
with the same parameters $c_{i,j}, \alpha_{i,j} \in [0, 1]$, 
$\underset{n \rightarrow \infty}{\lim}  
|F(\bar{p}_{1, n}, \ldots, \bar{p}_{m, n}) - F(\bar{q}_{1, n}, \ldots, \bar{q}_{m, n})| = 0$.
\end{enumerate}
}\end{defin}

\begin{exam}\label{the length functions are ct-continuous} {\rm
It is easy to see that for every $\beta \in (0, 1)$ the aggregation function 
$\mr{length}^{-\beta} : [0, 1]^{<\omega} \to [0, 1]$ defined by $\mr{length}^{-\beta}(\bar{p}) = |\bar{a}|^{-\beta}$
is ct-continuous for all parameters $c_1, \ldots, c_m$, $\alpha_1, \ldots, \alpha_m \in [0, 1]$ and all $m \in \mbbN^+$.
Also, it is not hard to prove that aggregation function `tsum' (or ``truncated sum'') used in
Example~\ref{example of page rank} is ct-continuous for all parameters
$c_1, \ldots, c_m$, $\alpha_1, \ldots, \alpha_m \in [0, 1]$ and all $m \in \mbbN^+$
}\end{exam}

\noindent
The proof of Proposition~6.3 in \cite{KW1} proves the next result because in the case of
am and gm the argument works even if some $\alpha_j$ is zero.

\begin{prop}\label{aggregation functions which are ct-continuous} \cite{KW1}
(i) For all $m \in \mbbN^+$ and all $c_1, \ldots, c_m, \alpha_1, \ldots, \alpha_m \in [0, 1]$
the aggregation functions am and gm are ct-continuous with respect to the parameters 
$c_1, \ldots, c_m$, $\alpha_1, \ldots, \alpha_m$.\\
(ii) For all $m \in \mbbN^+$, all $c_1, \ldots, c_m \in [0, 1]$ and all $\alpha_1, \ldots, \alpha_m \in (0, 1]$
(so $\alpha_i > 0$ for all $i$) the aggregation functions max and min are ct-continuous with respect to the
parameters $c_1, \ldots, c_m, \alpha_1, \ldots, \alpha_m$.
\end{prop}

\subsection{Uniform point continuity}

In this section we define our second notion of continuity, up-continuity, and for this we need
to associate every $\bar{p} \in [0, 1]^{<\omega}$ with a function $f_{\bar{p}} : [0, 1] \to [0, 1]$.
Actually we do this in two ways, an ``ordered'' and an ``unordered'' way. 
Common of both ways is that if $c \in \rng(\bar{p})$ and the proportion of coordinates in $\bar{p}$ that are equal to $c$
is $\alpha$, then $c$ belongs to the range/image of $f_{\bar{p}}$ and $f^{-1}(c)$ is a union of intervals 
such that the sum of the lengths of the intervals is $\alpha$.

Then we show how ct-continuity and up-continuity are related and give examples.

\begin{defin}\label{definition of associated function}{\bf (Functional representations of sequences)} {\rm
Let $n \in \mbbN^+$ and let $\bar{p} = (p_1, \ldots, p_n) \in [0, 1]^n$. 
To $\bar{p}$ we associate two functions from $[0, 1]$ to $[0, 1]$.
\begin{enumerate}
\item Define $\mff_{\bar{p}}$, which we call the {\em ordered functional representation of $\bar{p}$}, as follows: 
For every $a \in [0, 1/n)$, let $\mff_{\bar{p}}(a) = p_1$, 
for every $i = 1, \ldots, n-1$ and every $a \in [i/n, (i+1)/n)$, let $f(a) = p_{i+1}$ and finally let $f(1) = p_n$.

\item Define $\mfg_{\bar{p}}$, which we call the {\em unordered functional representation of $\bar{p}$}, as follows:
Let $\bar{p}' = (p'_1, \ldots, p'_n)$ be a reordering of $\bar{p}$ such that, for all $i = 1, \ldots, n-1$, 
$p'_i \leq p'_{i+1}$ and let $\mfg_{\bar{p}} = \mff_{\bar{p}'}$. 
\end{enumerate}
}\end{defin}

\begin{defin}\label{definition of the metric}{\bf (Pseudometrics on $\big([0, 1]^{<\omega}\big)^k$)} 
{\rm
\begin{enumerate}
\item First we recall the $L_1$ and $L_\infty$ norms: for every (bounded and integrable) $f : [0, 1] \to \mbbR$ they are defined as
\[
\| f \|_1 = \int_{[0,1]} | f(x) |dx \qquad \text{ and} \qquad
\| f \|_\infty = \sup\{|f(a)| : a \in [0, 1]\}.
\]
\item For $\bar{p}, \bar{q} \in [0, 1]^{<\omega}$ we define 
\begin{align*}
&\mu_1^u(\bar{p}, \bar{q}) =  \|\mfg_{p} - \mfg_{\bar{q}}\|_1, \\
&\mu_\infty^o(\bar{p}, \bar{q}) =  \|\mff_{\bar{p}} - \mff_{\bar{q}}\|_\infty.
\end{align*}

\item For arbitrary $k > 1$ we can define a functions on $\big([0, 1]^{<\omega}\big)^k$, also denoted 
$\mu_1^u$ and $\mu_\infty^o$
(to avoid making notation more complicated),
as follows: For all $(\bar{p}_1, \ldots, \bar{p}_k), (\bar{q}_1, \ldots, \bar{q}_k) \in \big([0, 1]^{<\omega}\big)^k$
let 
\begin{align*}
&\mu_1^u\big((\bar{p}_1, \ldots, \bar{p}_k), (\bar{q}_1, \ldots, \bar{q}_k)\big) = 
\max\big(\mu_1^u(\bar{p}_1, \bar{q}_1), \ldots, \mu_1^u(\bar{p}_k, \bar{q}_k)\big), \text{ and}\\
&\mu_\infty^o\big((\bar{p}_1, \ldots, \bar{p}_k), (\bar{q}_1, \ldots, \bar{q}_k)\big) = 
\max\big(\mu_\infty^u(\bar{p}_1, \bar{q}_1), \ldots, \mu_\infty^u(\bar{p}_k, \bar{q}_k)\big).
\end{align*}
\end{enumerate}
}\end{defin}

\noindent
It follows that $\mu_1^u$ and  $\mu_\infty^o$ are symmetric and satisfy the triangle inequality so they are pseudometrics on $[0, 1]^{<\omega}$.
None of them is a metric since it can happen that $\mu_1^u(\bar{p}, \bar{q}) = 0$ 
and $\bar{p} \neq \bar{q}$.
For example, if $\bar{p} = (0, 1/2, 1)$ and $\bar{q} = (0, 0, 1/2, 1/2, 1, 1)$
then $\mu_1^u(\bar{p}, \bar{q}) = 0$.
Also observe that for all $\bar{p}, \bar{q} \in [0, 1]^{<\omega}$, 
$\mu_1^u(\bar{p}, \bar{q}), \mu_\infty^o(\bar{p}, \bar{q}) \leq 1$.

\begin{defin}\label{definition of uniformly continuous aggregation function on a set} 
{\bf (Asymptotic uniform continuity on a set)} 
{\rm
Let $F : \big([0, 1]^{<\omega}\big)^k \to [0, 1]$ be an aggregation function and let 
$\mu$ be any of the the pseudometrics defined in Definition~\ref{definition of the metric}.
Also let $X \subseteq \big([0, 1]^{<\omega}\big)^k$.
We say that $F$ is {\em asymptotically uniformly continuous on $X$}
if for every $\varepsilon > 0$ there are $n$ and $\delta > 0$ such that 
if $(\bar{p}_1, \ldots, \bar{p}_k), (\bar{q}_1, \ldots, \bar{q}_k) \in X$,
$|\bar{p}_i|, |\bar{q}_i| \geq n$ for all $i$ and 
$\mu_1^u\big((\bar{p}_1, \ldots, \bar{p}_k), (\bar{q}_1, \ldots, \bar{q}_k)\big) < \delta$, then 
$\big|F(\bar{p}_1, \ldots, \bar{p}_k) - F(\bar{q}_1, \ldots, \bar{q}_k)\big| < \varepsilon$.
}\end{defin}

\begin{defin}\label{definition of up-continuous} 
{\bf (Uniform point continuity)} 
{\rm
An aggregation function $F : \big([0, 1]^{<\omega}\big)^m \to [0, 1]$ is called 
{\em up-continuous (uniformly point continuous) with respect to the parameters}
 $c_{i, j}, \alpha_{i, j} \in [0, 1]$, $i = 1, \ldots, m$ and $j = 1, \ldots, k_i$,
if the following two conditions hold:
\begin{enumerate}
\item For all sufficiently small $\delta > 0$,
$F$ is asymptotically uniformly continuous on $X_1 \times \ldots \times X_m$ where, for each $i = 1, \ldots, m$,
\begin{align*}
X_i = &\big\{\bar{p} \in [0, 1]^{<\omega} : \rng(\bar{p}) \subseteq \{c_{i,1}, \ldots, c_{i,k_i}\} 
\text{ and, for each $j = 1, \ldots, k_i$,} \\
&\text{there are between $(\alpha_{i,j} - \delta)|\bar{p}|$ and $(\alpha_{i,j} + \delta)|\bar{p}|$ coordinates in $\bar{p}$ } \\
&\text{which equals $c_{i,j}$} \big\}.
\end{align*}

\item For all $\varepsilon > 0$ there are $\delta > 0$ and $n_0$ such that if, for $i = 1, \ldots, m$,
$\bar{p}_i, \bar{q}_i \in [0, 1]^{<\omega}$ and
\begin{enumerate}
\item $|\bar{p}_i| = |\bar{q}_i| > n_0$, 
\item $\mu_\infty^o(\bar{p}_i, \bar{q}_i) < \delta$, 
\item $\rng(\bar{p}_i) \subseteq \{c_{i,1}, \ldots, c_{i,k_i}\}$, and
\item for each $j = 1, \ldots, k_i$, there are between $(\alpha_{i,j} - \delta)|\bar{p}_i|$ and 
$(\alpha_{i,j} + \delta)|\bar{p}_i|$ coordinates in $\bar{p}_i$ which equal $c_{i, j}$,
\end{enumerate}
then $|F(\bar{p}_1, \ldots, \bar{p}_m) - F(\bar{q}_1, \ldots, \bar{q}_m)| < \varepsilon$.
\end{enumerate}
}\end{defin}

\begin{exam}\label{example that mu-u is up-continuous} {\rm
It is straightforward to verify that, for all $m, k_1, k_2 \in \mbbN^+$, and all $c_{i, j}, \alpha_{i, j} \in [0, 1]$
for $i = 1, 2$, $j = 1, \ldots, k_i$, 
the pseudometric $\mu_1^u : \big([0, 1]^{<\omega}\big)^2 \to [0, 1]$
is up-continuous with respect to the parameters $c_{i, j}, \alpha_{i, j} \in [0, 1]$.
}\end{exam}

\noindent
The next proposition relates ct-continuity and up-continuity.

\begin{prop}\label{proposition relating the notions of continuity}
Let $F : \big([0, 1]^{<\omega}\big)^k \to [0, 1]$ be an aggregation function.\\
(i) If, for some $m, k_1, \ldots, k_m \in \mbbN^+$,
$F$ is up-continuous with respect to the parameters $c_{i, j}, \alpha_{i, j} \in [0, 1]$,
$i = 1, \ldots, m$, $j = 1, \ldots, k_i$, 
then $F$ is $ct$-continuous with respect to the same parameters $c_{i, j}, \alpha_{i, j} \in [0, 1]$.\\
(ii) Let $m, k_1, \ldots, k_m \in \mbbN^+$ and $c_{i, j}, \alpha_{i, j} \in [0, 1]$ for
$i = 1, \ldots, m$ and $j = 1, \ldots, k_i$.
If there is $\delta > 0$ such that, for all $\alpha'_{i, j} \in (\alpha_{i, j} - \delta, \alpha_{i, j} + \delta)$,
$i = 1, \ldots, m$ and $j = 1, \ldots, k_i$,
$F$ is ct-continuous with respect to the parameters $c_{i, j}, \alpha'_{i, j}$,
$i = 1, \ldots, m$ and $j = 1, \ldots, k_i$, then $F$ is up-continuous with respect to the 
parameters $c_{i, j}, \alpha_{i, j}$, $i = 1, \ldots, m$ and $j = 1, \ldots, k_i$.
\end{prop}

\noindent
Observe that we immediately get the following:

\begin{cor}\label{corollary to proposition relation the notions of continuity}
If $F : \big([0, 1]^{<\omega}\big)^m \to [0, 1]$ is an aggregation function then
the following are equivalent:
\begin{enumerate}
\item[(i)] $F$ is ct-continuous with respect to every choice of parameters $c_{i, j}, \alpha_{i, j} \in [0, 1]$,
$i = 1, \ldots, m$, $j = 1, \ldots, k_i$.
\item[(ii)] $F$ is up-continuous with respect to every choice of parameters $c_{i, j}, \alpha_{i, j} \in [0, 1]$,
$i = 1, \ldots, m$, $j = 1, \ldots, k_i$.
\end{enumerate}
\end{cor}

\noindent
{\bf Proof of Proposition~\ref{proposition relating the notions of continuity}.}
In order to avoid very cluttered notation and make the ideas of the proof more evident we only
prove the proposition in the case when $F$ is unary, i.e. when $F :[0, 1]^{<\omega} \to [0, 1]$.
The general case is proved in the same way except that we need to 
``book keep'' more sequences and parameters.

(i) Suppose that $F$ is up-continuous with respect to the parameters $c_1, \ldots, c_k$, $\alpha_1, \ldots, \alpha_k$.
Let $\bar{p}_n \in [0 1]^{<\omega}$ and $\bar{q}_n \in [0, 1]^{<\omega}$, 
$n \in \mbbN$, be convergence testing sequences for the parameters $c_1, \ldots, c_k$ and $\alpha_1, \ldots, \alpha_k$.
We need to show that $\lim_{n\to\infty}|F(\bar{p}_n) - F(\bar{q}_n)| = 0$.
Let $I_1, \ldots, I_k$ be disjoint open intervals such that $c_j \in I_j$ for $j = 1, \ldots, k$.
Let $\bar{p}_{n, i}$ denote the $i$th coordinate of $\bar{p}_n$ and similarly for $\bar{q}_{n, i}$.
As $\bar{p}_n$ and $\bar{q}_n$ are convergence testing sequences for
$c_1, \ldots, c_k$ and $\alpha_1, \ldots, \alpha_k$ it follows that there is $n_0$ such that if 
$n \geq n_0$, then $\rng(\bar{p}_n), \rng(\bar{q}_n) \subseteq I_1 \cup \ldots \cup I_k$, and
\begin{equation}\label{properties of the two convergence testing sequences}
\underset{n\to\infty}{\lim} \frac{|\{ i : p_{n, i} \in I_j\}|}{|\bar{p}_n|} \ = \
\underset{n\to\infty}{\lim} \frac{|\{ i : q_{n, i} \in I_j\}|}{|\bar{q}_n|} \  = \ \alpha_j.
\end{equation}
Since we are only considering the limit, we can assume without loss of generality that $n_0 =1$.  
Define the sequences  $\bar{p}'_n$ and $\bar{q}'_n$ by setting $p'_{n, i} =  c_j$ 
if  $p_{n, i} \in I_j$ (recall that different $I_j$ are disjoint),
and $q'_{n, i} =  c_j$ if $q_{n, i} \in I_j$. 

By condition~(1) of Definition~\ref{definition of up-continuous}
of up-continuity, for every $\varepsilon > 0$ there are $\delta > 0$ and $n_0$ such that if
$n > n_0$ and $\mu_1^u(\bar{p}'_n, \bar{q}'_n) < \delta$, then $|F(\bar{p}'_n) - F(\bar{q}'_n)| < \varepsilon$.
From~(\ref{properties of the two convergence testing sequences}) and the definition of $\bar{p}'_n$ and $\bar{q}'_n$
it follows that for every $\delta > 0$ we have $\mu_1^u((\bar{p}'_n, \bar{q}'_n) < \delta$ for all sufficiently large $n$.
Hence $\lim_{n\to\infty}|F(\bar{p}'_n) - F(\bar{q}'_n)| = 0$.

It now suffices to prove that 
\[
\underset{n \rightarrow \infty}{\lim}  |F(\bar{p}_n) - F(\bar{p}'_n)| \ = \ 
\underset{n \rightarrow \infty}{\lim}  |F(\bar{q}_n) - F(\bar{q}'_n)| = 0.
\]
We only show that the first limit is 0, since the the second limit is treated in the same way.
We will see that this is a consequence of Condition~(2) of 
Definition~\ref{definition of up-continuous} 
of up-continuity.  
By the choice of $\bar{p}_n$ (as convergence testing) and construction of $\bar{p}'_n$ we have
$|\bar{p}_n| = |\bar{p}'_n|$, $\lim_{n\to\infty}|\bar{p}_n| = \infty$, 
$\rng(\bar{p}'_n) \subseteq \{c_1, \ldots, c_k\}$, and for every $\delta > 0$ there is $n_0$ such that,
$\mu_\infty^o(\bar{p}_n, \bar{p}'_n) < \delta$ if $n > n_0$ 
(because $I_1, \ldots, I_k$ can be chosen with diameter at most $\delta$).
From~(\ref{properties of the two convergence testing sequences}) 
and the construction of $\bar{p}'_n$ it follows that for every $\delta > 0$ there is $n_0$ such that for all
$n > n_0$ and $j = 1, \ldots, k$, the number of coordinates in $\bar{p}'_n$ which equals $c_j$ is between
$(\alpha_j - \delta)|\bar{p}'_n|$ and $(\alpha_j + \delta)|\bar{p}'_n|$.
As $F$ is up-continuous with respect to $c_1, \ldots, c_k$ and $\alpha_1, \ldots, \alpha_k$,
it follows from condition~(2) of 
Definition~\ref{definition of up-continuous} 
that $\lim_{n\to\infty}  |F(\bar{p}_n) - F(\bar{p}'_n)| = 0$.

(ii) We prove the contrapositive. Suppose that $F$ is not up-continuous with respect to
the parameters $c_j, \alpha_j$, $j = 1, \ldots, k$.
We will show that there are arbitrarily small $\delta > 0$ and 
$\alpha'_j  \in (\alpha_j - \delta, \alpha_j + \delta)$ such that 
$F$ is not ct-continuous with respect to the parameters
$c_j, \alpha'_j$, $j = 1, \ldots, k$.
As $F$ is not up-continuous with respect to the mentioned parameters
either condition~(1) or condition~(2) of 
Definition~\ref{definition of up-continuous}
of up-continuity fails for the same parameters.

First suppose that condition~(1) of Definition~\ref{definition of up-continuous} fails
for the parameters $c_j, \alpha_j$, $j = 1, \ldots, k$.
Then there are arbitrarily small $\delta > 0$ such that $F$ is not asymptotically uniformly
continuous on 
\begin{align*}
X^\delta = &\big\{\bar{p} \in [0, 1]^{<\omega} : \rng(\bar{p}) \subseteq \{c_1, \ldots, c_k\} 
\text{ and, for each $j = 1, \ldots, k$,} \\
&\text{there are between $(\alpha_j - \delta)|\bar{p}|$ and $(\alpha_j + \delta)|\bar{p}|$ 
coordinates in $\bar{p}$ } \\
&\text{which equals $c_j$} \big\}.
\end{align*}
Let $\delta > 0$.
Without loss of generality we may assume that $\delta$ is small enough so that
if $\alpha_j > 0$ then $\alpha_j - \delta > 0$, and if $\alpha_j < 1$ then $\alpha_j + \delta < 1$.
As $F$ is not asymptotically uniformly continuous on $X^\delta$, there
is $\varepsilon > 0$ such that for all $m, N \in \mbbN^+$ there are 
$\bar{p}_{m, N}, \bar{q}_{m, N} \in X^\delta$ such that $\mu_1^u(\bar{p}_{m, N}, \bar{q}_{m, N}) < 1/m$,
$|\bar{p}_{m, N}|, |\bar{q}_{m, N}| > N$ and $|F(\bar{p}_{m, N}) - F(\bar{q}_{m, N})| > \varepsilon$.
We now define two convergence testing sequences the parameters of which will become clear later.

Define $\bar{p}'_1 = \bar{p}_{1, 1}$, $\bar{q}'_1 = \bar{q}_{1, 1}$, and for $n > 1$,
$\bar{p}'_n = \bar{p}_{n, N}$ and $\bar{q}'_n = \bar{q}_{n, N}$ where $N$ is
larger than $|\bar{p}'_m|$ and $|\bar{q}'_m|$ for all $m < n$.
Then $|\bar{p}'_n| < |\bar{p}'_{n+1}|$, $|\bar{q}'_n| < |\bar{q}'_{n+1}|$,
$\mu_1^u(\bar{p}'_n, \bar{q}'_n) < 1/n$ and $|F(\bar{p}'_n) - F(\bar{q}'_n)| > \varepsilon$
for all $n$.
Let $\bar{p}'_{n, i}$ denote the $i$th coordinate of $\bar{p}'$.
Since $[0, 1]^k$ is a compact topological space, the sequence of $k$-tuples
\begin{equation}\label{the sequence of proportions}
\big( |\{i : \bar{p}'_{n, i} = c_1\}|/|\bar{p}'_n|, \ldots,  |\{i : \bar{p}'_{n, i} = c_k\}|/|\bar{p}'_n| \big)
\end{equation}
has a convergent subsequence with limit $(\alpha'_1, \ldots, \alpha'_k)$, say.
Without loss of generality we may assume that the convergent subsequence is the whole sequence,
so~(\ref{the sequence of proportions}) converges to $(\alpha'_1, \ldots, \alpha'_k)$ as $n \to \infty$.
By the definition of $X^\delta$ and the construction of $\bar{p}'_n$ and $\bar{q}'_n$ we have
$\alpha_j - \delta \leq \alpha'_j \leq \alpha_j + \delta$ for all $j = 1, \ldots, k$.

We now show that both the sequence $\bar{p}'_n$ and the sequence $\bar{q}'_n$ is convergence
testing for the parameters $c_1, \ldots, c_k$ and $\alpha'_1, \ldots, \alpha'_k$.
We have already observed that $|\bar{p}'_n| < |\bar{p}'_{n+1}|$, $|\bar{q}'_n| < |\bar{q}'_{n+1}|$
so condition~(1) in Definition~\ref{definition of convergence testing} 
of convergence testing sequence is satisfied for both sequences.
Since $\rng(\bar{p}'_n), \rng(\bar{q}'_n) \subseteq \{c_1, \ldots, c_k\}$ (because both sequences belong to $X^\delta$)
it follows that for all open intervals $I_1, \ldots, I_k$ of $[0, 1]$
(with respect to the induced topology on $[0, 1]$) such that $c_j \in I_j$ for $j = 1, \ldots, k$,
$\rng(\bar{p}'_n), \rng(\bar{q}'_n) \subseteq \bigcup_{j=1}^k I_j$.
By the choice of $\alpha'_j$, $j = 1, \ldots, k$, it follows that 
\begin{equation}\label{the sequence of proportions again}
\lim\limits_{n \rightarrow \infty} \frac{\left| \{ i \leq |\bar{p}'_n| : p_{n,i} \in I_j \} \right| }{|\bar{p}'_n|} = \alpha_j.
\end{equation}
Hence the sequence $\bar{p}'_n$ is convergence testing for parameters $c_1, \ldots, c_k$ and $\alpha'_1, \ldots, \alpha'_k$.
Since $\mu_1^u(\bar{p}'_n, \bar{q}'_n) < 1/n$ it follows that~(\ref{the sequence of proportions again})
holds if $\bar{p}'_n$ is replaced by $\bar{q}'_n$, so $\bar{q}'_n$ is convergence testing for the same parameters.
As $|F(\bar{p}'_n) - F(\bar{q}'_n)| > \varepsilon$ for all $n$ we conclude that $F$ is not ct-continuous
with respect to the sequence of parameters $c_1, \ldots, c_k, \alpha_1, \ldots, \alpha_k$.

Now suppose that condition~(2) in 
Definition~\ref{definition of up-continuous}
of up-continuity fails for the sequence of parameters $c_1, \ldots, c_k, \alpha_1, \ldots, \alpha_k$.
Then there is $\varepsilon > 0$ such that for all $n \in \mbbN^+$ there are 
$\bar{p}_n, \bar{q}_n \in [0, 1]^{<\omega}$ such that 
\begin{enumerate}
\item[(a)] $|\bar{p}_n| = |\bar{q}_n| > n$, 
\item[(b)] $\mu_\infty^o(\bar{p}_n, \bar{q}_n) < 1/n$, 
\item[(c)] $\rng(\bar{p}_n) \subseteq \{c_1, \ldots, c_k\}$, and
\item[(d)] for each $j = 1, \ldots, k$, there are between $(\alpha_j - \delta)|\bar{p}_n|$,
$(\alpha_j + \delta)|\bar{p}_n|$ coordinates in $\bar{p}_n$ which equal $c_j$, and 
\item[(e)] $|F(\bar{p}_n) - F(\bar{q}_n)| \geq \varepsilon$.
\end{enumerate}
Without loss of generality we can also assume that $|\bar{p}_n| < |\bar{p}_{n+1}|$ and $|\bar{q}_n| < |\bar{q}_{n+1}|$
for all $n$.
From~(c) and~(d) it follows that $\bar{p}_n$, $n \in \mbbN^+$, is a convergence testing sequence for the 
parameters $c_1, \ldots, c_k$ and $\alpha_1, \ldots, \alpha_k$.
From~(b) it follows that for every $j = 1, \ldots, k$ and every open interval $I_j$ around $c_j$
there is $N$ such that if $n > N$ then an entry of $\bar{p}_n$ lies in $I_j$ if and only if the corresponding entry
of $\bar{q}_n$ lies in $I_j$. Hence $\bar{q}_n$, $n \in \mbbN^+$, is also a convergence testing sequence
for the parameters $c_1, \ldots, c_k$ and  $\alpha_1, \ldots, \alpha_k$.
From~(e) we now conclude that $F$ is not ct-continuous with respect to the parameters
$c_1, \ldots, c_k$ and  $\alpha_1, \ldots, \alpha_k$.
\hfill $\square$

\begin{exam}\label{am and gm are continuous for all parameters} {\rm
From Proposition~\ref{aggregation functions which are ct-continuous} 
and Proposition~\ref{proposition relating the notions of continuity} it follows that
the aggregation functions
`am' (arithmetic mean) and `gm' (geometric mean) are up-continuous with respect to all parameters, and that 
the aggregation functions `max' and `min'
are, for all $m, k_1, \ldots, k_m \in \mbbN^+$,
all $c_{i, j} \in [0, 1]$ and all $\alpha_{i, j} \in (0, 1]$,
up-continuous with respect to the parameters $c_{i, j}$ and $\alpha_{i, j}$.
From earlier examples and 
Proposition~\ref{proposition relating the notions of continuity}
it follows that `$\mr{length}^{-\beta}$' (where $\beta \in (0, 1)$) 
and `tsum' (``truncated sum'') are up-continuous with respect to all 
parameters.
More examples of aggregation functions that are ct-continuous and up-continuous with respect to
all parameters, or possibly with the requirement that $\alpha_{i, j} > 0$ for all $i$ and $j$,
are found in \cite{KW1} and \cite{KW2}.
}\end{exam}

\noindent
In the next two examples we use the following notation:
If $\alpha \in [0, 1]$ and $\bar{p} \in [0, 1]^{<\omega}$, then 
$\#(\alpha, \bar{p})$ is the number of coordinates in $\bar{p}$ which are equal to $\alpha$.
The remaining two examples consider aggregation functions obtained from 
generalized quantifiers, as discussed in 
Section~\ref{Expressing quantifiers with aggregation functions},
and show that rather mild conditions on a sequence of parameters imply that 
such aggregation functions are ct-continuous, respectively
up-continuous, with respect to that sequence of parameters.

\begin{exam}\label{example of proportional quantifiers} {\rm
In \cite{KL} Keisler and Lotfallah proved results about almost sure elimination of probability quantifiers
and convergence in probability.
Let $\beta \in (0, 1)$. The probability quantifier ``the proportion of $x$ satisfying $\ldots$ is at least $\beta$'' corresponds to
the following generalized quantifier
(in the sense of Definition~\ref{definition of generalized quantifier}):
\[
Q = \big\{ (D, X) : \ D \neq \es \text{ and } |X|/|D| \geq \beta \big\}.
\]
As we saw in the proof of 
Proposition~\ref{generalized quantifiers can be expressed in PLA*},
$Q$ can be represented by the aggregation function $F : [0, 1]^{<\omega} \to [0, 1]$ defined as follows:
\[
F(\bar{p}) = 
\begin{cases}
1 \qquad \text{ if } \#(1, \bar{p})/|\bar{p}| \geq \beta, \\
0 \qquad \text{ otherwise.}
\end{cases}
\]
Suppose that $\beta > 0$ and that $c_1, \ldots, c_k, \alpha_1, \ldots, \alpha_k \in  [0, 1]$.
It is straightforward to verify that $F$ is ct-continuous with respect to the parameters
$c_1, \ldots, c_k, \alpha_1, \ldots, \alpha_k$
if either
\begin{enumerate}
\item $c_1, \ldots, c_k < 1$, or if 

\item there are $m, i_1, \ldots, i_m \in \mbbN^+$ such that $c_{i_1}, \ldots, c_{i_m}$ is an enumeration
of all $c_i$ such that $c_i = 1$ and $\alpha_{i_1} + \ldots + \alpha_{i_m} \neq \beta$.
\end{enumerate}
Proposition~\ref{proposition relating the notions of continuity}
implies that under the same conditions $F$ is up-continuous with respect to the parameters
$c_1, \ldots, c_k, \alpha_1, \ldots, \alpha_k$.
}\end{exam}

\begin{exam}\label{example of rescher quantifier} {\rm
Consider the Rescher quantifier:
\begin{align*}
R = \big\{ (D, X_1, X_2) : \ |X_1| \leq |X_2| \big\}.
\end{align*}
It can be represented by the aggregation function
$G :  \big([0, 1]^{<\omega}\big)^2 \to [0, 1]$ defined as follows:
\begin{align*}
G(\bar{p}, \bar{q}) = 
\begin{cases}
1 \qquad \text{ if } \#(1, \bar{p}) \leq \#(1, \bar{q}), \\
0 \qquad \text{ otherwise.}
\end{cases}
\end{align*}
For $i = 1, 2$ and $j = 1, \ldots, k_i$ let $c_{i, j}, \alpha_{i, j} \in [0, 1]$.
If $c_{i, j} < 1$ for all $i$ and $j$ then it is easy to verify that $G$ is ct-continuous with respect to $c_{i, j}$ and $\alpha_{i, j}$.
Now suppose that $c_{1, j} = 1$ for some $j$ or that $c_{2, j} = 1$ for some $j$.
Let $j_1, \ldots, j_s$ enumerate all $j$ such that $c_{1, j} = 1$ and let
$l_1, \ldots, l_t$ enumerate all $l$ such that $c_{2, l} = 1$.
It is straightforward to verify that if 
\[
\alpha_{1, j_1} + \ldots + \alpha_{1, j_s} \neq \alpha_{2, l_1} + \ldots + \alpha_{2, l_t}
\]
then $G$ is ct-continuous with respect to $c_{i, j}$ and $\alpha_{i, j}$. 
(We consider the left side sum to be 0 if no $c_{1, j}$ exists which is 1, and we consider the
right side sum to be 0 if no $c_{2, j}$ exists which is 1.)
Proposition~\ref{proposition relating the notions of continuity}
implies that under the same conditions $G$ is up-continuous with respect to the parameters $c_{i, j}$ and $\alpha_{i, j}$.

The reader may verify that ct-continuity and up-continuity of the aggregation function which represents
the Härtig quantifier $H = \big\{ (D, X_1, X_2) : \ |X_1| = |X_2| \big\}$ can be characterized in the same way.
}\end{exam}

\section{Asymptotic equivalence to basic formulas}\label{Asymptotic equivalence to basic formulas}

\noindent
In this section we state and prove our main result, Theorem~\ref{main result}.
We begin by fixing some assumptions and notation for the rest of the section and then discuss the assumptions and
results of this section.

\begin{notation} {\rm
For each $n \in \mbbN^+$, $D_n$ is a finite set, also called a {\em domain}, and 
$|D_{n+1}| > |D_n|$ for all $n$.
By $\sigma$ we denote a finite and relational signature and $\mbW_n$ denotes a set of
(not necessarily all)
$\sigma$-structures with domain $D_n$.
By $\mbbP_n$ we denote a probability distribution on $\mbW_n$.
}\end{notation}

\noindent
A commonly considered context is when $D_n = \{1, \ldots, n\}$ and $\mbW_n$ is the set of all $\sigma$-structures
with domain $D_n$. However we want to allow for more generality.
For example one may be interested in $\sigma$-structures with domain $D_n$ such that all these structures have some
common properties. For example, we may restrict some relation symbols in $\sigma$ to be interpreted as relations
with degree\footnote{
A relation has degree at most $d$ if the Gaifman graph (also called the primal graph) corresponding to the relation
has no vertex with degree more than $d$.
}
at most 20 (say), or we could restrict a binary relation symbol $R \in \sigma$ to be interpreted as
a tree, a partial order, an equivalence relation, perhaps with some additional properties. 
We can, but need not, require that, for some $\sigma' \subset \sigma$ every structure in $\mbW_n$ has the same reduct to $\sigma'$.
If such restrictions are imposed on the $\sigma$-structures considered with domain $D_n$, then the sequence of
cardinalities $|D_n|$, $n \in \mbbN^+$, may not contain all (sufficiently large) positive integers.
Note also that we do not require that $D_n \subseteq D_{n+1}$.

With this quite general set-up our goal is to isolate properties 
(given by Assumption~\ref{assumptions on the basic logic})
of ``simpler'' 0/1-valued sublogics 
$L_0, L_1 \subseteq PLA^*(\sigma)$ such that if $L_0$ and $L_1$ have these properties, 
then a formula $\varphi(\bar{x}) \in PLA^*(\sigma)$ is asymptotically equivalent 
(see Definition~\ref{definition of asymptotic equivalence})
to an ``$L_0$-basic formula'', provided that every aggregation function
in $\varphi$ is up-continuous (or ct-continuous) with respect to some parameters which are
determined only by $(\mbbP_n : n \in \mbbN^+)$ and the subformulas of $\varphi$.
By an ``$L_0$-basic formula'' we mean a formula of the form $\bigwedge_{i=1}^k (\varphi_i(\bar{x}) \to c_i)$
where $\varphi_i \in L_0$ and $c_i \in [0, 1]$ for all $i$
and $\lim_{n\to\infty} \mbbP_n\big(\forall \bar{x} \bigvee_{i=1}^k \varphi_i(\bar{x})\big) = 1$.
If some $\psi(\bar{x}) \in PLA^*(\sigma)$ is asymptotically equivalent to an $L_0$-basic formula
$\varphi(\bar{x}) = \bigwedge_{i=1}^k (\varphi_i(\bar{x}) \to c_i)$,
then, for sufficiently large $n$, $\mcA \in \mbW_n$ and $\bar{a} \in (D_n)^{|\bar{x}|}$, 
the value $\mcA(\psi(\bar{a}))$
is, with high probability, close to the value $\mcA(\varphi(\bar{a}))$
which is determined by the values $\mcA(\varphi_i(\bar{a}))$, $i = 1, \ldots, k$, where $\varphi_i \in L_0$
are 0/1-valued.
In earlier work like \cite{Jae98a, KW1, KW2} and \cite{KL, Kop20} (if the later are reformulated in the context of $PLA^*$)
$L_0$ can be taken to be the set of complete consistent conjunctions of first-order literals. 
But this choice of $L_0$ is not always possible. 

For 0/1-valued logics such as first-order logic 
the notion of {\em almost sure/everywhere equivalence} between two formulas $\varphi(\bar{x})$ and $\psi(\bar{x})$
means that the probability that the value of $\forall \bar{x}(\varphi(\bar{x}) \leftrightarrow \psi(\bar{x}))$ is 1
tends to 1 as $n\to\infty$.
Since we consider a logic with (truth) values in the unit interval $[0, 1]$ we need to generalize 
the notion of almost sure equivalence and we do it as follows:

\begin{defin}\label{definition of asymptotic equivalence} {\bf (Equivalence and asymptotic equivalence)} {\rm
Let $\varphi(\bar{x}), \psi(\bar{x}) \in PLA^*(\sigma)$.\\
(i) We say that $\varphi(\bar{x})$ and $\psi(\bar{x})$ are {\em equivalent} if for
for every finite $\sigma$-structure $\mcA$ and every $\bar{a} \in A^{|\bar{x}|}$, 
$\mcA\big(\varphi(\bar{a})\big) = \mcA\big(\psi(\bar{a})\big)$.\\
(ii) We say that $\varphi(\bar{x})$ and $\psi(\bar{x})$ are 
{\em asymptotically equivalent (with respect to $(\mbbP_n : n \in \mbbN^+)$)}
if for every $\varepsilon > 0$ 
\[
\lim_{n\to\infty} \mbbP_n\Big(\Big\{\mcA \in \mbW_n : \text{ for all } \bar{a} \in (D_n)^{|\bar{x}|}, \
\big| \mcA\big(\varphi(\bar{a})\big) - \mcA\big(\psi(\bar{a})\big) \big| \leq \varepsilon \Big\}\Big) \ = \ 1.
\]
}\end{defin}

\noindent
Since we have fixed a sequence $\mbbP_n$, $n \in \mbbN^+$, of probability distributions for the rest of
the section we will just say that two formulas are asymptotically equivalent, omitting 
``with respect to $(\mbbP_n, n \in \mbbN^+)$''.

\begin{notation} {\rm
If $\varphi \in PLA^*(\sigma)$ is a formula without free variables then
\[
\mbbP_n\big(\varphi\big) = \mbbP_n\big(\big\{ \mcA \in \mbW_n : \mcA(\varphi) = 1 \big\}\big).
\]
}\end{notation}

\begin{defin}\label{L-basic subformula} {\bf (Special kinds of formulas)} {\rm
(i) A formula in $PLA^*(\sigma)$ such that no aggregation function occurs in it is called {\em aggregation-free}.\\
(i) Let $L \subseteq PLA^*(\sigma)$. We say that $L$ is {\em $0/1$-valued} if for every formula
$\varphi(\bar{x}) \in L$, every finite $\sigma$-structure $\mcA$, and every $\bar{a} \in A^{|\bar{x}|}$,
$\mcA(\varphi(\bar{a}))$ is either 0 or 1.\\
(ii) Let $L \subseteq PLA^*$ and suppose that $L$ is $0/1$-valued. 
A formula of $PLA^*$ is called {\em $L$-basic (formula)} if it has the form
$\bigwedge_{i=1}^k \big(\varphi_i(\bar{x}) \to c_i\big)$ where $\varphi_i \in L$ and $c_i \in [0, 1]$ for all $i = 1, \ldots, k$,
and $\lim_{n\to\infty}\mbbP_n\big(\forall \bar{x} \bigvee_{i=1}^k \varphi_i(\bar{x})\big) = 1$.
}\end{defin}

\noindent
When speaking of $L$-basic formulas we typically think of $L$ as a ``simple'' sublogic of $PLA^*$.
For example, $L$ could be the set of all consistent conjunctions of first-order literals.
Suppose for the moment that $L \subseteq L' \subseteq PLA^*(\sigma)$
and $L$ is the set of all consistent conjunctions of first-order literals.
If we can then prove that
(with respect to $\mbbP_n, n \in\mbbN^+$) every $\varphi(\bar{x}) \in L'$ is asymptotically equivalent to
an $L$-basic formula, then it means that, for large $n$ and with high probability, the value $\mcA(\varphi(\bar{a}))$
is close to a number which is determined by the first-order literals that are satisfied by $\bar{a}$ in $\mcA$.

The intuition behind the next assumption is that for some sets of formulas $L_0, L_1 \subseteq PLA^*(\sigma)$ 
and all $\varphi(\bar{x}, \bar{y}) \in L_0$ there is a 
set $L_{\varphi(\bar{x}, \bar{y})} \subseteq L_1$ of formulas defining some ``allowed'' conditions,
and if $\chi(\bar{x}, \bar{y}) \in L_{\varphi(\bar{x}, \bar{y})}$, 
then the fraction $|\varphi(\bar{a}, \mcA) \cap \chi(\bar{a}, \mcA)| / |\chi(\bar{a}, \mcA)|$ is with high probability close to a 
number $\alpha$ that depends only on the involved formulas and the sequence $\mbbP_n$.
The assumptions in the main results of for example \cite{Jae98a, Kop20, KW1, KW2} imply that 
Assumption~\ref{assumptions on the basic logic} below holds if $L_0$ is the set of all formulas which are 
(consistent) conjunctions
of first-order literals (and more generally quantifier-free first-order formulas). 
In the same articles $L_1$ can, depending on the result, be either $\{\top\}$
(where, for any $\bar{x}$ and $\bar{y}$, we can view $\top$ as a formula
$\chi(\bar{x}, \bar{y})$ which has value 1 for any $\bar{a}$ and $\bar{b}$
from any structure), the set of all
complete and consistent conjunctions of identities ($x = y$) and nonidentities ($x \neq y$), or the set of all
consistent conjuctions of identities and nonidentities;
in the same results $L_{\varphi(\bar{x}, \bar{y})} = L_1$ for every $\varphi(\bar{x}, \bar{y}) \in L_0$.

\begin{assump}\label{assumptions on the basic logic}{\rm
Suppose that $L_0 \subseteq PLA^*(\sigma)$ and $L_1 \subseteq PLA^*(\sigma)$ 
are $0/1$-valued and that the following conditions hold:
\begin{enumerate}
\item For every aggregation-free $\varphi(\bar{x}) \in PLA^*(\sigma)$ there is 
an $L_0$-basic formula $\varphi'(\bar{x})$ such that
$\varphi$ and $\varphi'$ are asymptotically equivalent.

\item For every $m \in \mbbN^+$ and all
$\varphi_1(\bar{x}, \bar{y}), \ldots, \varphi_m(\bar{x}, \bar{y}) \in L_0$, 
there are $L_{\varphi_j(\bar{x}, \bar{y})} \subseteq L_1$ for $j = 1, \ldots, m$
such that if $\chi_j(\bar{x}, \bar{y}) \in L_{\varphi_j(\bar{x}, \bar{y})}$
for $j = 1, \ldots, m$, then 
there are $s, t \in \mbbN^+$, $\theta_i(\bar{x}) \in L_0$, $\alpha_{i, j} \in [0, 1]$,
for $i = 1, \ldots, s$, $j = 1, \ldots, m$,
and $\chi'_i(\bar{x}) \in L_0$, for $i = 1, \ldots, t$,
such that for every $\varepsilon > 0$ and $n$ there is $\mbY^\varepsilon_n \subseteq \mbW_n$ such that
$\lim_{n\to\infty}\mbbP_n(\mbY^\varepsilon_n) = 1$ and 
for every $\mcA \in \mbY^\varepsilon_n$ the following conditions hold:
\begin{align*}
&(a) \ \mcA \models \forall \bar{x} \bigvee_{i=1}^s \theta_i(\bar{x}), \\
&(b) \ \text{if $i \neq j$ then } 
\mcA \models \forall \bar{x} \neg(\theta_i(\bar{x}) \wedge \theta_j(\bar{x})), \\
&(c) \ \mcA \models \forall \bar{x}
\Big(\Big(\bigvee_{i=1}^m \neg\exists\bar{y}\chi_i(\bar{x}, \bar{y})\Big) \leftrightarrow
\Big(\bigvee_{i=1}^t \chi'_i(\bar{x})\Big)\Big), \text{ and}\\
&(d) \ \text{for all $i = 1, \ldots, s$ and $j = 1, \ldots, m$, if $\bar{a} \in (D_n)^{|\bar{x}|}$, and
$\mcA \models \theta_i(\bar{a})$,} \\
&\text{ then } (\alpha_{i, j} - \varepsilon)|\chi_j(\bar{a}, \mcA)| \ \leq \
|\varphi_j(\bar{a}, \mcA) \cap \chi_j(\bar{a}, \mcA)|
\ \leq \ (\alpha_{i, j} + \varepsilon)|\chi_j(\bar{a}, \mcA)|.
\end{align*}
\end{enumerate}
}\end{assump}

\noindent
From now on we assume that Assumption~\ref{assumptions on the basic logic} holds,
that is, we have fixed two $0/1$-valued sublogics $L_0, L_1 \subseteq PLA^*(\sigma)$ such that 
conditions~(1) and~(2) of  Assumption~\ref{assumptions on the basic logic} are satisfied.

\begin{lem}\label{L-0 basic formulas are closed under connectives}
Suppose that $\msfC : [0, 1]^k \to [0, 1]$ is a continuous connective.
If $\varphi_1(\bar{x}),$ $\ldots,$ $\varphi_k(\bar{x}) \in$ $PLA^*(\sigma)$ and, for $i = 1, \ldots, k$,
$\varphi_i(\bar{x})$ is asymptotically equivalent to an $L_0$-basic formula $\psi_i(\bar{x})$,
then $\msfC(\varphi_1(\bar{x}), \ldots, \varphi_k(\bar{x}))$ is asymptotically equivalent
to an $L_0$-basic formula.
\end{lem}

\noindent
{\bf Proof.}
Suppose that the assumptions of the lemma are satisfied.
Since $\msfC$ is continuous and, for $i = 1, \ldots, k$, $\varphi_i(\bar{x})$ and $\psi_i(\bar{x})$
are asymptotically equivalent it follows that 
$\msfC(\varphi_1(\bar{x}), \ldots, \varphi_k(\bar{x}))$ and
$\msfC(\psi_1(\bar{x}), \ldots, \psi_k(\bar{x}))$
are asymptotically equivalent.
Since the formula $\msfC(\psi_1(\bar{x}), \ldots, \psi_k(\bar{x}))$ is aggregation-free
it follows, by part~(1) of 
Assumption~\ref{assumptions on the basic logic}, 
that it is equivalent to an $L_0$-basic formula $\psi'(\bar{x})$.
It now follows that $\msfC(\varphi_1(\bar{x}), \ldots, \varphi_k(\bar{x}))$ and
$\psi'(\bar{x})$ are asymptotically equivalent.
\hfill $\square$

\begin{defin}\label{y-frequency parameters} {\bf (Frequency parameters of an $L_0$-basic formula)} {\rm
Let $L_0$ and $L_1$ be sublogics of $PLA^*$ that satisfy 
Assumption~\ref{assumptions on the basic logic}.
Let $\psi(\bar{x}, \bar{y})$ be an $L_0$-basic formula, so $\psi(\bar{x}, \bar{y})$ has the form
$\bigwedge_{i=1}^s\big(\psi_i(\bar{x}, \bar{y}) \to c_i\big)$
where $\psi_i \in L_0$ and $c_i \in [0, 1]$ for $i = 1, \ldots, s$, and
$\lim_{n\to\infty}\mbbP_n\big(\forall \bar{x}, \bar{y} \bigvee_{i=1}^s \psi_i(\bar{x}, \bar{y})\big) = 1$.
Suppose that $\chi(\bar{x}, \bar{y}) \in \bigcap_{i = 1}^s L_{\psi_i(\bar{x}, \bar{y})}$.

It follows from 
Assumption~\ref{assumptions on the basic logic} 
that there are
$\theta_1(\bar{x}), \ldots, \theta_t(\bar{x}) \in L_0$ 
and $\alpha_{i, j} \in [0, 1]$, for $i = 1, \ldots, t$ and $j = 1, \ldots, s$,
such that 
for every $\varepsilon > 0$ there is $\mbY^\varepsilon_n \subseteq \mbW_n$
such that $\lim_{n\to\infty}\mbbP_n\big(\mbY^\delta_n\big) = 1$ and for all $\mcA \in \mbY^\varepsilon_n$
\begin{align*}
&\mcA \models \forall \bar{x} \bigvee_{i=1}^t\theta_i(\bar{x}), \\
&\text{if $i \neq j$ then } 
\mcA \models \forall \bar{x} \neg(\theta_i(\bar{x}) \wedge \theta_j(\bar{x})),
\end{align*}
and, for all $i = 1, \ldots, t$ and all $j = 1, \ldots, s$, if 
$\bar{a} \in (D_n)^{|\bar{x}|}$,
$\mcA \models \theta_i(\bar{a})$ and  $\chi(\bar{a}, \mcA) \neq \es$, then
\[
\alpha_{i, j} - \varepsilon \leq \frac{|\psi_j(\bar{a}, \mcA) \cap \chi(\bar{a}, \mcA)|}{|\chi(\bar{a}, \mcA)|}
\leq \alpha_{i, j} + \varepsilon.
\]
Fix some $i \in \{1, \ldots, t\}$.
Let $c \in \{c_1, \ldots, c_s\}$ and let $j_1, \ldots, j_r$ enumerate, without repetition, all indices $i$ 
such that $c_i = c$. 
Let $\beta_c = \alpha_{i, j_1} + \ldots + \alpha_{i, j_r}$. 
Let $d_1, \ldots, d_{s'}$ be an enumeration of $\{c_1, \ldots, c_s\}$ without repetition.
Then we call $(d_1, \ldots, d_{s'}, \beta_{d_1}, \ldots, \beta_{d_{s'}})$ the {\em sequence of
$\bar{y}$-frequency parameters of $\psi$ relative to $\chi$ and $\theta_i$}
(and it is unique up to the order of the $d_i$ and the corresponding $\beta_{d_i}$).
A {\em sequence of $\bar{y}$-frequency parameters of $\psi$ relative to $\chi$} is, by definition,
a sequence of $\bar{y}$-frequency parameters of $\psi$ relative to $\chi$ and $\theta_i$, for some $i$.

The above definition can be generalized, for any $r > 1$, to $r$ $L_0$-basic formulas
$\psi_1(\bar{x}, \bar{y}), \ldots, \psi_r(\bar{x}, \bar{y})$.
Suppose that, for each $i = 1, \ldots, r$, $\psi_i$ has the form 
\[
\bigwedge_{k=1}^{s_i} \big(\psi_{i, k}(\bar{x}, \bar{y}) \to c_{i, k}\big),
\]
so each $\psi_{i, k}$ belongs to $L_0$.
Suppose that, for $i = 1, \ldots, r$, 
$\chi_i(\bar{x}, \bar{y}) \in \bigcap_{k=1}^{s_i} L_{\psi_{i, k}(\bar{x}, \bar{y})}$.
By Assumption~\ref{assumptions on the basic logic} 
there are 
$\theta_1(\bar{x}), \ldots, \theta_t(\bar{x}) \in L_0$ 
and $\alpha_{i, j, k} \in [0, 1]$ 
for $i = 1, \ldots, t$, $j = 1, \ldots, r$ and $k = 1, \ldots, s_j$, such that 
for every $\varepsilon > 0$ there is $\mbY^\varepsilon_n \subseteq \mbW_n$
such that $\lim_{n\to\infty}\mbbP_n\big(\mbY^\delta_n\big) = 1$ and for all $\mcA \in \mbY^\varepsilon_n$
\begin{align*}
&\mcA \models \forall \bar{x} \bigvee_{i=1}^t\theta_i(\bar{x}), \\
&\text{if $i \neq j$ then } 
\mcA \models \forall \bar{x} \neg(\theta_i(\bar{x}) \wedge \theta_j(\bar{x})),
\end{align*}
and, for all $i = 1, \ldots, t$ and all $j = 1, \ldots, r$, if 
$\bar{a} \in (D_n)^{|\bar{x}|}$,
$\mcA \models \theta_i(\bar{a})$ and  $\chi_j(\bar{a}, \mcA) \neq \es$, then
\[
\alpha_{i, j, k} - \varepsilon \leq \frac{|\psi_{j, k}(\bar{a}, \mcA) \cap \chi_j(\bar{a}, \mcA)|}{|\chi_j(\bar{a}, \mcA)|}
\leq \alpha_{i, j, k} + \varepsilon.
\]
{\em The $r$-tuple of sequences of $\bar{y}$-frequency parameters of $(\psi_1, \ldots, \psi_r)$
relative to $(\chi_1, \ldots, \chi_r)$ and $\theta_i$} is, by definition,
$(\bar{p}_1, \ldots, \bar{p}_r)$, where $\bar{p}_j$ is the sequence of
$\bar{y}$-frequency parameters of $\psi_j$ relative to $\chi_j$ and $\theta_i$.
An $r$-tuple is an {\em $r$-tuple of $\bar{y}$-frequency parameters of $(\psi_1, \ldots, \psi_r)$
relative to $(\chi_1, \ldots, \chi_r)$} if, for some $i$, it is the sequence of
$\bar{y}$-frequency parameters of $(\psi_1, \ldots, \psi_r)$
relative to $(\chi_1, \ldots, \chi_r)$ and $\theta_i$.
}\end{defin}

\noindent
The next lemma follows more or less directly from 
Definition~\ref{y-frequency parameters}
and states the properties of frequency parameters
that will be relevant later.

\begin{lem}\label{proportion of entries equal to a number in a sequence}
Let $\psi(\bar{x}, \bar{y})$ denote the $L_0$-basic formula
$\bigwedge_{i=1}^s \big(\psi_i(\bar{x}, \bar{y}) \to c_i\big)$ and suppose that
$\chi(\bar{x}, \bar{y}) \in \bigcap_{i=1}^s L_{\psi_i(\bar{x}, \bar{y})}$.
Also let $\theta_1(\bar{x}), \ldots, \theta_t(\bar{x}) \in L_0$ be formulas as in
Definition~\ref{y-frequency parameters},
and let $(d_{i, 1}, \ldots, d_{i, s_i}, \beta_{i, 1}, \ldots, \beta_{i, s_i})$ be the sequence of 
$\bar{y}$-frequency parameters of $\psi$ relative to $\chi$ and $\theta_i$. Then:
\begin{enumerate}
\item For each $i = 1, \ldots, t$, $d_{i, 1}, \ldots, d_{i, s_i}$ is an enumeration of 
$\{c_1, \ldots, c_s\}$ without repetition.

\item Suppose that $\mbY^\varepsilon_n$ is like in
Definition~\ref{y-frequency parameters} and $\mcA \in \mbY^\varepsilon_n$.
Let $i \in \{1, \ldots, t\}$ and suppose that $\bar{a} \in (D_n)^{|\bar{x}|}$, $\chi(\bar{a}, \mcA) \neq \es$,
$\mcA \models \theta_i(\bar{a})$, and let 
\[
\bar{p} = \big(\mcA\big(\psi(\bar{a}, \bar{b})\big) : 
\bar{b} \in (D_n)^{|\bar{x}|} \text{ and } \mcA \models \chi(\bar{a}, \bar{b})\big).
\]
Then
\begin{enumerate}
\item $\rng(\bar{p}) \subseteq \{c_1, \ldots, c_s\} = \{d_{i, 1}, \ldots, d_{i, s_i}\}$.

\item For all $j = 1, \ldots, s_i$, then number of coordinates in $\bar{p}$ which are equal to $d_{i, j}$
is between $(\beta_{i, j} - s \varepsilon)|\chi(\bar{a}, \mcA)|$ and 
$(\beta_{i, j} + s \varepsilon)|\chi(\bar{a}, \mcA)|$.
\end{enumerate}
\end{enumerate}
\end{lem}

\begin{theor}\label{main result}
Let $L_0$ and $L_1$ be sublogics of $PLA^*(\sigma)$ that satisfy 
Assumption~\ref{assumptions on the basic logic}.
Let $F : \big([0, 1]^{<\omega}\big)^m \to [0, 1]$ be an $m$-ary aggregation function, 
let $\varphi_i(\bar{x}, \bar{y}) \in PLA^*(\sigma)$, for $i = 1, \ldots, m$, and suppose that each 
$\varphi_i(\bar{x}, \bar{y})$
is asymptotically equivalent to an $L_0$-basic formula $\psi_i(\bar{x}, \bar{y})$ where
$\psi_i$ has the form $\bigwedge_{k=1}^{s_i} (\psi_{i, k}(\bar{x}, \bar{y}) \to c_{i, k})$.
Suppose that for $i = 1, \ldots, m$, 
$\chi_i(\bar{x}, \bar{y}) \in \bigcap_{k=1}^{s_i} L_{\psi_{i, k}(\bar{x}, \bar{y})}$.
Let $\varphi(\bar{x})$ denote the $PLA^*(\sigma)$-formula
\[
F\big(\varphi_1(\bar{x}, \bar{y}), \ldots, \varphi_m(\bar{x}, \bar{y}) : \bar{y} : 
\chi_1(\bar{x}, \bar{y}), \ldots, \chi_m(\bar{x}, \bar{y})\big).
\]
(i) If $F$ is up-continuous with respect to every $m$-tuple of $\bar{y}$-frequency parameters of
$(\psi_1, \ldots, \psi_m)$ relative to $(\chi_1, \ldots, \chi_m)$, then $\varphi(\bar{x})$
is asymptotically equivalent to an $L_0$-basic formula.\\
(ii) Suppose that there is $\delta > 0$ such that $F$ is ct-continuous with respect to every
$m$-tuple of (sequences of) parameters $(\bar{p}_1, \ldots, \bar{p}_m)$ such that the following holds:
\begin{enumerate}
\item[] For each $i = 1, \ldots, m$, if
$\bar{p}_i = (c_{i, 1}, \ldots, c_{i, k_i}, \alpha_{i, 1}, \ldots, \alpha_{i, k_i})$, then there are 
$\beta_{i, j} \in [0, 1]$ such that 
$\alpha_{i, j} \in (\beta_{i, j} - \delta, \beta_{i, j} + \delta)$ for $j = 1, \ldots,k_i$, and
if 
\[
\bar{q}_i =  (c_{i, 1}, \ldots, c_{i, k_i}, \beta_{i, 1}, \ldots, \beta_{i, k_i}),
\]
then $(\bar{q}_1, \ldots, \bar{q}_m)$ is an $m$-tuple of $\bar{y}$-frequency parameters of 
$(\psi_1, \ldots, \psi_m)$ relative to $(\chi_1, \ldots, \chi_m)$.
\end{enumerate}
Then $\varphi(\bar{x}, \bar{y})$ is asymptotically equivalent to an $L_0$-basic formula.
\end{theor}

\begin{exam}\label{concrete example of the main result}{\rm
As a concrete example of an application of Theorem~\ref{main result} let us consider the following context.
Let $\sigma$ be a finite relational signature and let
$\mbW_n$ be the set of all $\sigma$-structures with domain $\{1, \ldots, n\}$.
Suppose that $\mbbP_n$ is the uniform probability distribution on $\mbW_n$ for all $n$.
Let $L_0$ be the set of all conjunctions of first-order literals (over $\sigma$) and let
$L_1 = L_0$.
In this setting, an $L_0$-basic formula is one of the form $\bigwedge_{i=1}^s (\theta_i(\bar{x}) \to c_i)$
where $c_i \in [0, 1]$ and $\theta_i(\bar{x}) \in L_0$ for all $i = 1, \ldots, s$, 
and $\forall \bar{x} \bigvee_{i=1}^s \theta_i(\bar{x})$ is a valid sentence.

One can easily prove, as in \cite{KW1}, that every aggregation-free formula is equivalent to an $L_0$-basic formula
so part~(1) Assumption~\ref{assumptions on the basic logic} holds.
From the proofs of the classical 0-1 law (in particular the arguments in \cite{Gleb})
it follows that also part~(2) of Assumption~\ref{assumptions on the basic logic} holds
and we can, with the notation of that assumption, let 
$L_{\varphi_j(\bar{x}, \bar{y})} = L_1$ for all 
$\varphi_j(\bar{x}, \bar{y}) \in L_0$.
It now follows from
Example~\ref{am and gm are continuous for all parameters}
and Theorem~\ref{main result} that if
$\varphi(\bar{x}, \bar{y})$ is an $L_0$-basic formula and $\chi(\bar{x}, \bar{y}) \in L_1$,
then the formulas
\[
\mr{am}\big( \varphi(\bar{x}, \bar{y}) : \bar{y} : \chi(\bar{x}, \bar{y}) \big)
\ \ \text{ and } \ \ 
\mr{gm}\big( \varphi(\bar{x}, \bar{y}) : \bar{y} : \chi(\bar{x}, \bar{y}) \big)
\]
(where am and gm are the arithmetic, respectively gemetric,  means) 
are asymptotically equivalent to $L_0$-basic formulas.
If we require that the formula called $\chi(\bar{x}, \bar{y})$ above is a
conjuction of formulas of the form $x = y$ or $x \neq y$, then,
by the results in \cite{KW1, KW2}, we get the same conclusion
for more probability distributions 
(than the uniform one) and for more aggregation functions (than am and gm).
In \cite{Kop26, KT}, which use Theorem~\ref{main result},
other choices of $L_0$ and $L_1$ turn out to be appropriate 
for the results aimed at there.\footnote{
The articles \cite{Kop26, KT} were written after this article but published before.}
}\end{exam}

\begin{rem}\label{applications of the main result} {\bf (Applications of Theorem~\ref{main result})} {\rm
(i) Suppose that Assumption~\ref{assumptions on the basic logic}
holds and $\varphi(\bar{x}) \in PLA^*(\sigma)$.
Part~(1) of Assumption~\ref{assumptions on the basic logic},
Lemma~\ref{L-0 basic formulas are closed under connectives}, and 
Theorem~\ref{main result},
tell that $\varphi(\bar{x})$ can be ``reduced'', step by step, to an asymptotically equivalent $L_0$-basic formula
as long as, in a step when an aggregation function $F$ is to be ``asymptotically eliminated'' we do not
encounter parameters with respect to which $F$ is not up-continuous.
The proofs of
Lemma~\ref{L-0 basic formulas are closed under connectives} and 
Theorem~\ref{main result}
indicate how the steps in such a reduction are carried out.

(ii) Now suppose that when trying to asymptotically eliminate an aggregation function $F$ as in
Theorem~\ref{main result} this $F$ is {\em not} up-continuous with respect to 
the parameters $c_{i, j}, \alpha_{i, j}$ that arise in this step. 
Then there may be another aggregation function $G$ which behaves ``almost'' like $F$ but
is up-continuous with respect to the parameters $c_{i, j}, \alpha_{i, j}$.
If $F$ is replaced by $G$ in the formula then the change of meaning of the new formula may be negligible in a particular context, 
but now the aggregation function $G$ can be eliminated.
For example, if $F$ is as in 
Example~\ref{example of proportional quantifiers}
and we encounter parameters $c_j, \alpha_j$ with respect to which $F$ is not up-continuous,
or equivalently not ct-continuous,
then we can consider an aggregation function $G$ which is defined just as $F$ except that we replace $\beta$
with some $\beta' \neq \beta$ which is very close to $\beta$.

(iii) As a special case of Theorem~\ref{main result} we have the following:
Suppose that Assumption~\ref{assumptions on the basic logic} holds for some $L_0, L_1 \subseteq PLA^*(\sigma)$
where $L_1 = \{\top\}$ and $L_{\varphi(\bar{x}, \bar{y})} = L_1$ for all $\varphi(\bar{x}, \bar{y}) \in L_0$.
Also let $coPLA^*(\sigma, L_0, L_1)$ be the set of all formulas $\varphi(\bar{x}) \in PLA^*(\sigma)$ such that 
every subformula of $\varphi(\bar{x})$ of the form 
\[
F\big(\varphi_1(\bar{x}, \bar{y}), \ldots, \varphi_m(\bar{x}, \bar{y}) : \bar{y} : 
\chi_1(\bar{x}, \bar{y}), \ldots, \chi_m(\bar{x}, \bar{y})\big).
\]
is such that $\chi_1, \ldots, \chi_m \in L_1$ (i.e. all $\chi_i$ are $\top$)
and $F$ is up-continuous with respect to all choices of parameters, or equivalently
(by Proposition~\ref{proposition relating the notions of continuity}),
ct-continuous with respect to all choices of parameters.
Examples of such aggregation functions are given by
Proposition~\ref{aggregation functions which are ct-continuous}
and Example~\ref{example that mu-u is up-continuous}.
Then
Lemma~\ref{L-0 basic formulas are closed under connectives}
and Theorem~\ref{main result} 
imply that {\em every formula in $coPLA^*(\sigma, L_0, L_1)$ is asymptotically equivalent
to an $L_0$-basic formula}. 
}\end{rem}

\begin{rem}\label{necessity of continuity} {\bf (Necessity of continuity)} {\rm
In Theorem~\ref{main result}
we cannot remove the assumption that $F$ is up-continuous, or ct-continuous, with respect to certain parameters.
This can be seen in different ways.
One way is the following, roughly explained.
The work in \cite{KL} coveres the context where $\mbW_n$ is the set of all $\sigma$-structures with domain
$D_n = \{1, \ldots, n\}$ and $\mbbP_n$ is the uniform probability distribution on $\mbW_n$.
The work in \cite{KL} shows that 
Assumption~\ref{assumptions on the basic logic} 
is satisfied if $L_0$ is the set of (consistent) conjunctions of first-order literals, $L_1 = \{\top\}$,
and $L_{\varphi(\bar{x}, \bar{y})} = L_1$ for all $\varphi(\bar{x}, \bar{y}) \in L_0$.
Let $PLA^*(\sigma, L_0, L_1)$ be the set of formulas $\varphi(\bar{x}) \in PLA^*(\sigma)$ such that
every subformula of $\varphi(\bar{x})$ of the form 
\[
F\big(\varphi_1(\bar{x}, \bar{y}), \ldots, \varphi_m(\bar{x}, \bar{y}) : \bar{y} : 
\chi_1(\bar{x}, \bar{y}), \ldots, \chi_m(\bar{x}, \bar{y})\big).
\]
is such that $\chi_1, \ldots, \chi_m \in L_1$ (i.e. all $\chi_i$ are $\top$).
As was explained in 
Example~\ref{example of proportional quantifiers}
every formula of the (0/1-valued) ``probability logic'' $L$ considered in \cite{KL} is equivalent
to a $PLA^*(\sigma, L_0, L_1)$-formula (if $\sigma$ contains all relation symbols used by $L$).
If Theorem~\ref{main result} would hold without the continuity assumption, 
then every $PLA^*(\sigma, L_0, L_1)$-formula
would be asymptotically equivalent to an $L_0$-basic formula.
It would follow (as $L_0$ contains only quantifier-free formulas) that every formula of the logic $L$ considered
in \cite{KL} is asymptotically equivalent to a quantifier-free first-order formula, and from this it would follow
that for every $\varphi \in L$, $\lim_{n\to\infty} \mbbP_n(\varphi)$ exists 
(see \cite[Corollary~4.10]{KL}).
But this is not the case, because as shown in \cite[Proposition~3.1]{KL} there is $\varphi \in L$ such that the limit
does not exist.
}\end{rem}

\subsection{Proof of Theorem~\ref{main result}}\label{proof of main result}

Part (ii) of Theorem~\ref{main result} is a consequence of part~(i) and 
Proposition~\ref{proposition relating the notions of continuity}.
For if the assumptions in part~(ii) hold, then $F$ is 
up-continuous with respect to every $m$-tuple of $\bar{y}$-frequency parameters of
$(\psi_1, \ldots, \psi_m)$ relative to $(\chi_1, \ldots, \chi_m)$, so by part~(i)
$\varphi(\bar{x})$ is asymptotically equivalent to an $L_0$-basic formula.

Hence it remains to prove part~(i) of 
Theorem~\ref{main result} and the rest of this section is devoted to this.
In order to make the arguments more clear by avoiding heavy notation we will assume that $m = 1$.
The general case is proved in the same way except that we need to keep track of more formulas and parameters.

Throughout the rest of this section we assume that $\psi(\bar{x}, \bar{y})$ denotes the $L_0$-basic formula
$\bigwedge_{i=1}^s \big(\psi_i(\bar{x}, \bar{y}) \to c_i\big)$,
so $\psi_i \in L_0$ for all $i$.
We also assume that $\chi(\bar{x}, \bar{y}) \in \bigcap_{i=1}^s L_{\psi_i(\bar{x}, \bar{y})}$.
By Assumption~\ref{assumptions on the basic logic} 
there are $\theta_1(\bar{x}), \ldots, \theta_t(\bar{x}) \in L_0$,
$\chi_1(\bar{x}), \ldots, \chi_k(\bar{x}) \in L_0$,
$\alpha_{i, j} \in [0, 1]$ for $i = 1, \ldots, t$ and $j = 1, \ldots, s$, and $\mbY^\delta_n \subseteq \mbW_n$
such that $\lim_{n\to\infty} \mbbP_n(\mbY^\delta_n) = 1$ and for all $\mcA \in \mbY^\delta_n$,
\begin{align*}
&(a) \ \mcA \models \forall \bar{x} \bigvee_{i=1}^t \theta_i(\bar{x}) = 1, \\
&(b) \ \text{if $i \neq j$ then } 
\mcA \models \forall \bar{x} \neg(\theta_i(\bar{x}) \wedge \theta_j(\bar{x})), \\
&(c) \ \mcA \models \forall \bar{x}
\Big(\neg\exists\bar{y}\chi(\bar{x}, \bar{y}) \leftrightarrow
\Big(\bigvee_{i=1}^k \chi_i(\bar{x})\Big)\Big), \text{ and}\\
&(d) \ \text{for all $i = 1, \ldots, t$,  if $\bar{a} \in (D_n)^{|\bar{x}|}$, and
$\mcA \models \theta_i(\bar{a})$,} \\
&\text{ then } (\alpha_{i, j} - \varepsilon)|\chi(\bar{a}, \mcA)| \ \leq \
|\psi_j(\bar{a}, \mcA) \cap \chi(\bar{a}, \mcA)|
\ \leq \ (\alpha_{i, j} + \varepsilon)|\chi(\bar{a}, \mcA)|.
\end{align*}
Furthermore, suppose that $F : [0, 1]^{<\omega} \to [0, 1]$ is an aggregation function that is 
up-continuous with respect to every sequence of $\bar{y}$-frequency parameters of $\psi$ relative to $\chi$.

Our first goal, achieved by 
Lemma~\ref{F(psi) asymptotic to basic formula} 
below, is to prove that $F\big(\psi(\bar{x}, \bar{y}) : \bar{y} : \chi(\bar{x}, \bar{y}) \big)$
is asymptotically equivalent to an $L_0$-basic formula.
Then we prove, as stated by Proposition~\ref{elimination of one aggregation function} below,
that if $\varphi(\bar{x}, \bar{y}) \in PLA^*(\sigma)$ and $\varphi(\bar{x}, \bar{y})$ is asymptotically equivalent to
$\psi(\bar{x}, \bar{y})$, then $F\big(\varphi(\bar{x}, \bar{y}) : \bar{y} : \chi(\bar{x}, \bar{y}) \big)$
is asymptotically equivalent to an $L_0$-basic formula.
This completes the proof of Theorem~\ref{main result}.

\begin{lem}\label{convergence relative to theta-i}
Fix an arbitrary index $1 \leq i \leq t$.
There is $d_i \in [0, 1]$, depending only on $\psi$, $\chi$, $\theta_i$, $F$, and $(\mbbP_n : n \in \mbbN^+)$ such that
for every $\varepsilon > 0$ there is $\delta > 0$ such that for all sufficiently large $n$, all $\mcA \in \mbY^\delta_n$, 
and all $\bar{a} \in (D_n)^{|\bar{x}|}$,
if $\chi(\bar{a}, \mcA) \neq \es$, and $\mcA \models \theta_i(\bar{a})$, then
\[
\big| \mcA\big(F\big( \psi(\bar{a}, \bar{y}) : \bar{y} : \chi(\bar{a}, \bar{y})\big)\big) - d_i \big| < \varepsilon.
\]
\end{lem}

\noindent
{\bf Proof.}
Fix $i \in \{1, \ldots, t\}$ and let $\varepsilon > 0$.
By the semantics of $PLA^*$ it suffices to show that 
there is $\delta > 0$ such that if $\mcA_1 \in \mbY^\delta_{n_1}$, $\mcA_2 \in \mbY^\delta_{n_2}$,
$\bar{a}_1 \in [n_1]^{|\bar{x}|}$, $\bar{a}_2 \in [n_2]^{|\bar{x}|}$,
$\chi(\bar{a}_1, \mcA_1) \neq \es$, $\chi(\bar{a}_2, \mcA_2) \neq \es$, 
$\mcA_1 \models \theta_i(\bar{a}_1)$, $\mcA_2 \models \theta_i(\bar{a}_2)$,
and, for $k = 1, 2$,
\[
\bar{p}_k = \big(\mcA_k(\psi(\bar{a}_k, \bar{b})) : \bar{b} \in [n_k]^{|\bar{y}|}
\text{ and } \mcA \models \chi(\bar{a}, \bar{b}) \big),
\]
then $|F(\bar{p}_1) - F(\bar{p}_2)| < \varepsilon$.

Let, according to 
Definition~\ref{y-frequency parameters},
$(d_{i, 1}, \ldots, d_{i, s_i}, \beta_{i, 1}, \ldots, \beta_{i, s_i})$ be the sequence of 
$\bar{y}$-frequency parameters of $\psi$ relative to $\chi$ and $\theta_i$.
By Lemma~\ref{proportion of entries equal to a number in a sequence},
$d_{i, 1}, \ldots, d_{i, s_i}$ is an enumeration of $\{c_1, \ldots, c_s\}$ without repetition,
for $k = 1, 2$,
$\rng(\bar{p}_k) \subseteq \{c_1, \ldots, c_s\}$,
and for all $j = 1, \ldots, s_i$ the number of coordinates in $\bar{p}_k$ which are equal to $d_{i, j}$ is
between $(\beta_{i, j} - s\delta)|\chi(\bar{a}, \mcA_k)|$
and $(\beta_{i, j} + s\delta)|\chi(\bar{a}, \mcA_k)|$.
It follows that $\mu_1^u(\bar{p}_1, \bar{p}_2) < C\delta$ where the constant $C$ depends only on $s_i$.

We assume that $F$ is up-continuous with respect to every sequence of $\bar{y}$-frequency parameters
of $\psi$ relative to $\chi$, so in particular $F$ is up-continuous
with respect to $(d_{i, 1}, \ldots, d_{i, s_i}, \beta_{i, 1}, \ldots, \beta_{i, s_i})$.
It follows from part~(1)
of the definition of up-continuity that if $\delta > 0$ is chosen small enough then
$|F(\bar{p}_1) - F(\bar{p}_2)| < \varepsilon$.
\hfill $\square$

\begin{lem}\label{F(psi) asymptotic to basic formula}
There is an $L_0$-basic formula $\psi'(\bar{x})$ such that for every $\varepsilon > 0$
there is $\delta > 0$ such that for all sufficiently large $n$, all $\mcA \in \mbY^\delta_n$, 
and all $\bar{a} \in (D_n)^{|\bar{x}|}$,
\[
\big| \mcA\big(F\big(\psi(\bar{a}, \bar{y}) : \bar{y} : \chi(\bar{a}, \bar{y})\big)\big) - 
\mcA\big(\psi'(\bar{a})\big) \big| < \varepsilon.
\]
\end{lem}

\noindent
{\bf Proof.}
Let $d_1, \ldots, d_t \in [0, 1]$ be as in 
Lemma~\ref{convergence relative to theta-i}.
By~(c) above, for all $\mcA \in \mbY^\delta_n$,
\[
\mcA \models
\forall \bar{x}\Big(\neg\exists\bar{y}\chi(\bar{x}, \bar{y}) \leftrightarrow
\big(\bigvee_{i=1}^k \chi_i(\bar{x})\big)\Big).
\]
Let $\psi'(\bar{x})$ be the formula
$\bigwedge_{i=1}^s (\theta_i(\bar{x}) \to d_i) \ \wedge \ 
\bigwedge_{j=1}^k (\chi_j(\bar{x}) \to 0)$,
so $\psi'$ is an $L_0$-basic formula.
We now verify that the claim of the lemma holds with this choice of $\psi'$.
First note that if $\bar{a} \in (D_n)^{|\bar{x}|}$, $\mcA \in \mbY^\delta_n$ and $\chi(\bar{a}, \mcA) = \es$,
then, for some $j$, $\mcA(\chi_j(\bar{a})) = 1$, so $\mcA(\chi_j(\bar{a}) \to 0) = 0$ for some $j$ and hence
\[
\mcA(\psi'(\bar{a})) = 0 = \mcA\big(F\big(\psi(\bar{a}, \bar{y}) : \bar{y} : \chi(\bar{a}, \bar{y})\big)\big).
\]
Now suppose that $\bar{a} \in (D_n)^{|\bar{x}|}$, $\mcA \in \mbY^\delta_n$ and $\chi(\bar{a}, \mcA) \neq \es$.
Hence $\mcA \not\models \chi_i(\bar{a})$ for all $i = 1, \ldots, k$ and thus
$\mcA\big(\bigwedge_{j=1}^k (\chi_j(\bar{x}) \to 0)\big) = 1$.
By~(a) and~(b) above there is a unique $i$ such that 
$\mcA \models \theta_i(\bar{a})$. 
Hence $\mcA\big(\theta_i(\bar{a}) \to d_i\big) = d_i$ and
$\mcA\big(\theta_{i'}(\bar{a}) \to d_{i'}\big) = 1$ for all $i' \neq i$.
It follows that 
$\mcA(\psi'(\bar{a})) = d_i$.
By Lemma~\ref{convergence relative to theta-i},
if $\delta > 0$ is sufficiently small we get 
\[
\big| \mcA\big(F\big( \psi(\bar{a}, \bar{y}) : \bar{y} : \chi(\bar{a}, \bar{y})\big)\big) - d_i \big| < \varepsilon.
\]
and the lemma now follows.
\hfill $\square$

\begin{prop}\label{elimination of one aggregation function}
Suppose that $\varphi(\bar{x}, \bar{y}) \in PLA^*(\sigma)$ and that $\varphi(\bar{x}, \bar{y})$
and $\psi(\bar{x}, \bar{y})$ are asymptotically equivalent.
Then 
$
F(\varphi(\bar{x}, \bar{y}) : \bar{y} : \chi(\bar{x}, \bar{y}))
$
is asymptotically equivalent to an $L_0$-basic formula.
\end{prop}

\noindent
{\bf Proof.}
Suppose that $\varphi(\bar{x}, \bar{y}) \in PLA^*(\sigma)$ and that $\varphi(\bar{x}, \bar{y})$
and $\psi(\bar{x}, \bar{y})$ are asymptotically equivalent.
By Lemma~\ref{F(psi) asymptotic to basic formula},
there is an $L_0$-basic formula $\psi'(\bar{x})$ such that, for all $\varepsilon > 0$, if 
$\delta > 0$ is small enough then for all sufficiently large $n$, all $\mcA \in \mbY^\delta_n$, 
and all $\bar{a} \in (D_n)^{|\bar{x}|}$,
\begin{equation}\label{F(psi) and psi' are close}
\big| \mcA\big(F\big(\psi(\bar{a}, \bar{y}) : \bar{y} : \chi(\bar{a}, \bar{y})\big)\big) - 
\mcA\big(\psi'(\bar{a})\big) \big| < \varepsilon/2.
\end{equation}
For any $\delta > 0$ define
\begin{align*}
\mbX^\delta_n = \big\{ \mcA \in \mbW_n : 
\text{ for all $\bar{a} \in (D_n)^{|\bar{x}|}$ and all $\bar{b} \in (D_n)^{|\bar{y}|}$},
\big|\mcA(\varphi(\bar{a}, \bar{b})) - \mcA(\psi(\bar{a}, \bar{b}))\big| < \delta \big\}.
\end{align*}
Since $\varphi(\bar{x}, \bar{y})$ and $\psi(\bar{x}, \bar{y})$ are asymptotically equivalent we have
$\lim_{n\to\infty}\mbbP_n(\mbX^\delta_n) =~1$.
Hence $\lim_{n\to\infty}\mbbP_n(\mbX^\delta_n \cap \mbY^\delta_n) = 1$.
Therefore it suffices to show that if $\delta > 0$ is small enough, $\mcA \in \mbX^\delta_n \cap \mbY^\delta_n$
and $\bar{a} \in (D_n)^{|\bar{x}|}$, then
\begin{equation}\label{to be proved in eliminating F in the general case}
\big| \mcA\big(F\big(\varphi(\bar{a}, \bar{y}) : \bar{y} : \chi(\bar{a}, \bar{y})\big)\big) - 
\mcA\big(\psi'(\bar{a})\big) \big| < \varepsilon.
\end{equation}
The inequality~(\ref{to be proved in eliminating F in the general case}) 
is a consequence of~(\ref{F(psi) and psi' are close})
and the following inequality:
\begin{equation}\label{difference between F(varphi) and F(psi)}
\big| \mcA\big(F\big(\varphi(\bar{a}, \bar{y}) : \bar{y} : \chi(\bar{a}, \bar{y})\big)\big) - 
\mcA\big(F\big(\psi(\bar{a}, \bar{y}) : \bar{y} : \chi(\bar{a}, \bar{y})\big)\big) < \varepsilon/2.
\end{equation}
So we need to prove that if $\delta > 0$ is sufficiently small and $n$ sufficiently large,
then~(\ref{difference between F(varphi) and F(psi)}) holds for
all $\mcA \in \mbX^\delta_n \cap \mbY^\delta_n$ and all $\bar{a} \in (D_n)^{|\bar{x}|}$.

Let $\mcA \in  \mbX^\delta_n \cap \mbY^\delta_n$ and $\bar{a} \in (D_n)^{|\bar{x}|}$.
First suppose that $\chi(\bar{a}, \mcA) = \es$. Then
\[
\mcA\big(F\big(\varphi(\bar{a}, \bar{y}) : \bar{y} :  \chi(\bar{a}, \bar{y})\big)\big) \ = \ 0 \ = \ 
\mcA\big(F\big(\psi(\bar{a}, \bar{y}) : \bar{y} : \chi(\bar{a}, \bar{y})\big)\big).
\]
Now suppose that $\chi(\bar{a}, \mcA) \neq \es$.
There is (by~(a) and~(b) above) a unique $i$ such that $\mcA \models \theta_i(\bar{a})$ and the following sequences are nonempty:
\begin{align*}
&\bar{p} = \big(\mcA\big(\varphi(\bar{a}, \bar{b})\big) : 
\bar{b} \in (D_n)^{|\bar{y}|} \text{ and } \mcA \models \chi(\bar{a}, \bar{b}) \big), \\
&\bar{q} = \big(\mcA\big(\psi(\bar{a}, \bar{b})\big) : 
\bar{b} \in (D_n)^{|\bar{y}|} \text{ and } \mcA \models \chi(\bar{a}, \bar{b}) \big).
\end{align*}
From 
Lemma~\ref{proportion of entries equal to a number in a sequence}
it follows that if $(c'_1, \ldots, c'_m, \alpha_1, \ldots, \alpha_m)$ is the sequence
of $\bar{y}$-frequency parameters of $\psi(\bar{x}, \bar{y})$ relative to $\chi(\bar{x}, \bar{y})$
and $\theta_i(\bar{x})$, then 
$\rng(\bar{q}) \subseteq \{c'_1, \ldots, c'_m\}$ and
the number of coordinates in $\bar{q}$ which are equal to $c'_j$ is between 
$(\alpha_j - s\delta)|\chi(\bar{a}, \mcA)|$ and
$(\alpha_j + s\delta)|\chi(\bar{a}, \mcA)|$.
Since $\mcA \in  \mbX^\delta_n$ it follows from the definition of  $\mu_\infty^o$ that
$\mu_\infty^o(\bar{p}, \bar{q}) < \delta$.

Since $F$ is up-continuous with respect to every sequence of 
$\bar{y}$-frequency parameters of $\psi(\bar{x}, \bar{y})$ relative to $\chi(\bar{x}, \bar{y})$,
it is, in particular, up-continuous with respect to the sequence 
$(c'_1, \ldots, c'_m, \alpha_1, \ldots, \alpha_m)$.
It now follows from condition~(2) of the definition of up-continuity with respect to the
sequence of parameters $(c'_1, \ldots, c'_m, \alpha_1, \ldots, \alpha_m)$ that 
if $\delta > 0$ is sufficiently small, 
then $|F(\bar{p}) - F(\bar{q})| < \varepsilon/2$.
This together with the definitions of $\bar{p}$ and $\bar{q}$ and the semantics of $PLA^*$ 
imply that~(\ref{difference between F(varphi) and F(psi)}) holds, so the proof is completed.
\hfill $\square$

\section{Conclusion}

\noindent
We have considered the quite expressive logic $PLA^*$ with truth values in the unit interval and 
which uses aggregation functions instead of
(generalized) quantifiers. 
We have shown that generalized quantifiers in the sense of Mostowski
\cite{Mos57, Mos} can be expressed by appropriate aggregation functions, so $PLA^*$ subsumes
(on finite structures)
first-order logic extendend by such generalized quantifiers.

We have also studied continuity properties of aggregation functions
and how the continuity properties of an aggregation function $F$ influence whether 
a $PLA^*$-formula that uses $F$ is asymptotically equivalent to a simpler formula
without occurences of $F$,
where `asymptotical equivalence' generalizes `almost sure equivalence' to the context of
continuous-valued logics.
The same approach has been used earlier in \cite{KW1, KW2}, but in more specific contexts.
This study has identified conditions which make sense in wider contexts than that of \cite{KW1, KW2} 
and the main result (Theorem~\ref{main result}) may therefore be useful in a variety of situations,
as already exemplified by \cite{Kop26, KT} which use the main result of this study.

\subsection*{Competing interests}

The authors have no competing interests as defined by Springer, or other interests that might be perceived to influence the results and/or discussion reported in this paper.

\subsection*{Research funding}

 This research did not receive funding.

\subsection*{Data availability}

This paper does not report data generation or analysis.

\subsection*{Author contributions}

Both authors have contributed to all results of the paper.

\end{document}